\documentclass{amsart}
\usepackage{amsmath}
\usepackage{comment}
\usepackage{amssymb}
\usepackage{amsthm}
\usepackage{dsfont}
\usepackage{lineno}
\usepackage{subfigure}
\usepackage{graphicx}
\usepackage{multirow}
\usepackage{color}
\usepackage{xspace}
\usepackage{pdfsync}
\usepackage{tikz-cd,mathtools}
\usepackage{hyperref} 
\usepackage{algorithmicx}
\usepackage{algpseudocode}
\usepackage{algorithm}
\usepackage{mathtools}
\usepackage{enumitem}
\graphicspath{{Figures/}}

\DeclareMathOperator*{\argmin}{argmin}

\DeclarePairedDelimiter\floor{\lfloor}{\rfloor}
\newcommand{\bq}{\begin{equation}}
\newcommand{\eq}{\end{equation}}
\newcommand{\R}{\mathbb{R}}

\newcommand{\abs}[1]{\left\vert#1\right\vert}

\newcommand{\G}{\mathcal{G}}
\newcommand{\bO}{\mathcal{O}}

\newcommand{\Dt}{\mathcal{D}}

\newcommand{\Sf}{\mathbb{S}^{2}}
\newcommand{\Tf}{\mathcal{T}}
\newcommand{\Zf}{\mathcal{Z}}

\algnewcommand{\LineComment}[1]{\State \(\triangleright\) #1}

\newtheorem{theorem}{Theorem}
\newtheorem{thm}{Theorem}

\theoremstyle{lemma}
\newtheorem{lemma}[theorem]{Lemma}

\newtheorem{definition}[theorem]{Definition}

\newtheorem{remark}[theorem]{Remark}

\newtheorem{hypothesis}[theorem]{Hypothesis}

\theoremstyle{remark}

\makeatletter
\newcommand\appendix@section[1]{%
\refstepcounter{section}%
\orig@section*{Appendix \@Alph\c@section: #1}%
}
\let\orig@section\section
\g@addto@macro\appendix{\let\section\appendix@section}
\makeatother

\begin{document}

\title[Adaptive Mesh Methods]{Adaptive Mesh Methods on Compact Manifolds via Optimal Transport and Optimal Information Transport}

\author{Axel G. R. Turnquist}
\address{Department of Mathematical Sciences, New Jersey Institute of Technology, University Heights, Newark, NJ 07102}
\email{agt6@njit.edu}

\thanks{The author was partially supported by an NSF GRFP 1849508 and by NSF DMS 1751996}

\begin{abstract}
Moving mesh methods are designed to redistribute a mesh in a regular way. This applied problem can be considered to overlap with the problem of finding a diffeomorphic mapping between density measures. In applications, an off-the-shelf grid needs to be restructured to have higher grid density in some regions than others. This should be done in a way that avoids tangling, hence, the attractiveness of diffeomorphic mapping techniques. For exact diffeomorphic mapping on the sphere a major tool used is Optimal Transport, which allows for diffeomorphic mapping between even non-continuous source and target densities. However, recently Optimal Information Transport was rigorously developed allowing for exact and inexact diffeomorphic mapping and the solving of a simpler partial differential equation.

In this manuscript, we perform the first side-by-side comparison of using Optimal Transport and Optimal Information Transport on the sphere for adaptive mesh problems. We introduce how to generalize these computations to more general manifolds. In this manuscript, we choose to perform this comparison with provably convergent solvers, which is generally challenging for either problem due to the lack of boundary conditions and lack of comparison principle in the partial differential equation formulation. It appears that Optimal Information Transport produces better results and is more easily generalizable for the moving mesh problem. In order to use Optimal Transport for moving mesh methods, further work on the accuracy of Optimal Transport solvers on the sphere and a framework for generalization to other manifolds must be done before this becomes a recommended method in challenging cases.
\end{abstract}

\date{\today}    
\maketitle

The diffeomorphic density matching problem has a long history in the imaging sciences. Classical methods consist of positive scalar image functions composed from the right by transformations~\cite{Younes}. That is, given a source probability density function $f_0$ and a target probability density function $f_1$, classical techniques compute $f_0 \circ T = f_1$. It is not possible in many cases to find diffeomorphic mappings. Furthermore, this formulation is not appropriate for moving mesh methods. Non-classical methods, like Optimal Transport and Optimal Information Transport, allow the transformation to act as a pushforward or pullback on the density, that is $f_0 = \left\vert D T \right\vert f_1 \circ T$ or $\left\vert D T \right\vert f_0 \circ T = f_1$, respectively. This generalization has particular benefit in that it allows for proving the existence of a diffeomorphic mapping over a much wider range of densities $f_0$ and $f_1$. Furthermore, it is clear that this is the formulation of diffeomorphic density matching that is appropriate for moving mesh methods, since it allows one to change the local density of mesh points.

Non-classical diffeomorphic density matching has found widespread applications in medical image registration  see, for example,~\cite{lungCT, HakerRegistration, chenozan, thoracic}, random sampling from a probability measure see, for example,~\cite{sampling, bayesian} and adaptive moving mesh methods, for example, in meteorology see the line of work in~\cite{budd1, budd2, BuddWilliams}, among other applications. These applications are necessarily more challenging from a theoretical perspective as well as a computational perspective when the geometry is non-Euclidean. In this manuscript, we perform the first comparison of Optimal Transport versus Optimal Information Transport for the moving mesh problem in the non-Euclidean case.

Optimal Transport density matching in non-Euclidean geometries is a field of study in which much progress on the theoretical side has been made recently, see the following for examples of this large line of work~\cite{mccann, LoeperReg, Loeper_OTonSphere, figalli, villanistability, loepervillani}.
Recently, the authors in~\cite{HT_OTonSphere, HT_OTonSphere2} introduced a numerical convergence framework and convergent numerical scheme for the PDE formulation of the optimal transport problem on the sphere which included, as one application, the traditional squared geodesic cost $c(x,y) = \frac{1}{2} d_{\mathbb{S}^2}(x,y)^2$. The convergence framework addressed situations where the solutions $u$ of the PDE formulation are known to be \textit{a priori} smooth (that is, at least $C^2$) as well as nonsmooth (but still $C^1$). This convergence framework required that the numerical schemes were consistent, stable in the average value, and monotone and achieved compactness via a Lipschitz constraint imposed on the PDE level. Preliminary results in~\cite{HT_OTonSphere2} showed the success of using Optimal Transport on the sphere for moving mesh methods. In this manuscript, we demonstrate the results in cases where the provably convergent numerical scheme for Optimal Transport on the sphere in~\cite{HT_OTonSphere2} is successful and unsuccessful for applications.

Optimal Information Transport was developed in the papers~\cite{modin, modin2}, which treat both exact and inexact diffeomorphic density matching based on a principal fiber bundle structure and the particularly nice way that horizontal geodesics in an information metric are pulled down to geodesics with respect to the Fisher-Rao metric. This allows for the exact diffeomorphic map between two smooth densities to be computed via solving a Poisson equation. Limitations to this approach are addressed as well, such as its apparent inappropriateness for image interpolation problems, by introducing computations with a non-compatible metric (which allows for the source density to vanish)~\cite{modin} in Euclidean space. Optimal Information Transport methods have subsequently found applications in image registration for medical applications, for example in~\cite{chenozan, thoracic} and also random sampling~\cite{sampling}, also in subsets of Euclidean space. In this manuscript, we use a provably convergent scheme for the Optimal Information Transport problem on the sphere and show its success for the moving mesh problem in challenging cases.

In this manuscript, we establish a class of provably convergent methods (with convergence rate guarantees) for computing the Optimal Information Transport on smooth, compact, and connected 2D manifolds $M$. We then provide, for the first time, a comparison of the adaptive mesh methods on these manifolds via Optimal Transport and Optimal Information Transport. In Section~\ref{sec:movingmesh} we introduce the applied problem of adaptive moving mesh methods. These have specific applications in the context of meteorology for the spherical geometry, but have applications also in medical image registration and computer graphics beyond the case of the sphere. In Section~\ref{sec:optimaltransport} we review briefly the problem of Optimal Transport on the sphere (and beyond) and highlight the regularity issues that arise which complicate the numerical convergence guarantees. In Section~\ref{sec:optimalinformationtransport} we outline how the Optimal Information Transport PDE is derived.s In Section~\ref{sec:numerics}, we show the grid setup and how consistent, monotone, and stable discrete operators are constructed for the PDE in Optimal Transport and Optimal Information Transport on the sphere. In Section~\ref{sec:implementation} we show computations for both problems on the sphere, highlighting the success of both in some cases and the failure of the Optimal Transport discretization in more challenging examples.

\textit{Acknowledgement}: I would like to express my gratitude to Brittany Froese Hamfeldt for fruitful discussions and for encouraging me to write this manuscript on my own.

\section{Moving Mesh Methods}\label{sec:movingmesh}
For this manuscript, our major application will be moving mesh methods, a type of adaptive mesh method that has found special usage in computational PDE techniques. The resulting computed mesh is then often used in high-resolution PDE computations, such as the Eady problem, which is a 2D vertical cross-section of an incompressible Bousinnesq fluid, which was treated in~\cite{budd1}. Fixing the number of grid points leads to simpler data structures (than some other adaptive mesh techniques) in the PDE solving step, as was stressed in the manuscripts~\cite{budd1, budd2, chacon}.

In the Euclidean setup, the moving mesh problem setup requires the ``physical" target domain $\Omega_p \subset \mathbb{R}^d$, where the PDE is posed, while the ``computational" input domain $\Omega_c \subset \mathbb{R}^d$ is usually chosen to be a uniform rectangular grid (in Euclidean space). In non-Euclidean geometries, one would like to use any simple off-the-shelf mesh generator for the computational mesh. Then, one would like to find a diffeomorphic mapping $T: \Omega_c \rightarrow \Omega_p$. Typically, in applications, the local density of grid points in the physical domain $\Omega_p$ is determined by, for example, the time scales involved in the solution of a fluid mechanics PDE. For example, in an evolving front or shock, it is desirable for the density of points in $\Omega_p$ to be greater than in the relatively unchanging parts of the solution, as in meteorological applications where it is desirable to have high resolution in areas of high precipitation~\cite{budd2}. This process could be done iteratively to get a very accurate solution of the PDE, by solving the shock on more and more resolved grids as one iterates. The way this is done is by feeding the information from the PDE (desired density) into a scalar monitor function $\mathcal{M}(y,t) > 0$ and then solving the change of variables formula:

\begin{equation}\label{eq:JacobianEquation}
\mathcal{M}(y,t) J(T(x)) = \theta(t)
\end{equation}
where $J$ is the Jacobian of the diffeomorphic mapping $T$ and $\int_{\Omega_p} \mathcal{M}(y,t) dy = \theta(t)$. Moving mesh methods require that the mapping be a diffeomorphism, which means that the Jacobian of the mapping $T$ satisfies $0<J(T)<\infty$ in the strong sense. This condition will prevent the mesh from tangling. An example of a monitor function constructed from information from a function $f(y,t)$ is the scaled arc-length function:

\begin{equation}
\mathcal{M}(y,t) =\sqrt{1 + \left\vert S \nabla_{y} f(y,t) \right\vert^2}
\end{equation}
where $S$ is a normalization factor, see~\cite{budd1}.

Applications of adaptive mesh techniques are also found in medical image registration and sampling techniques, which require computing diffeomorphic mapping between possibly very different density functions. Thus, we need to solve the problem with a broad variety of applications in mind.

The diffeomorphic density matching problem can be solved via the Jacobian equation, which is the change of variables. The existence of any such map is quite general, see~\cite{villani2}.

\begin{thm}[Jacobian]\label{thm:Jacobian}
Let $M$ be an $n$-dimensional Riemannian manifold. Let $\mu_0$ and $\mu_1$ be two probability measures on $M$, and let $T: M \rightarrow M$ be a measurable function such that $T_{\#} \mu_0 = \mu_1$. Let $\nu$ be a reference measure, of the form $\nu(dx) = e^{-V(x)} \text{vol}(dx)$, where $V$ is continuous and $\text{vol}$ is the volume measure on $M$. Further assume that $\mu_0(dx) = f_0 (x) \nu(dx)$ and $\mu_{1}(dy) = f_1 (y) \nu(dy)$, $T$ is injective, and the distributional derivative of the mapping $DT$ is a locally integrable function (this can be relaxed slightly). Then, $\mu_0$-almost surely,

\begin{equation}\label{eq:JacobianEquationVillani}
f_0(x) = f_1(T(x)) J(T(x))
\end{equation}
where $J$ is the Jacobian determinant of $T$ at $x$, defined weakly by:

\begin{equation}
J(T(x)) := \lim_{\epsilon \rightarrow 0} \frac{\nu \left[ T(B_{\epsilon}(x)) \right]}{\nu [B_{\epsilon}(x)]}
\end{equation}
\end{thm}

While the above theorem~\ref{thm:Jacobian} stipulates the conditions for which we have a distributional solution of the monitor equation~\eqref{eq:JacobianEquation}, it does not state that such a mapping $T$ is unique or when it is diffeomorphic. It does show that, however, for the same source and target mass distributions, there can be multiple solutions $T$. What we will find is that under fairly general assumptions on the source and target masses, we can find diffeomorphic mappings using the techniques of Optimal Transport and Optimal Information Transport, which both satisfy the Jacobian equation and give us diffeomorphic mappings, appropriate for our moving mesh methods.

The mapping arising from Optimal Transport is elected by further requiring that the map $T$ minimize the following integral:

\begin{equation}\label{eq:functional}
I = \int_{\Omega_c} \left\vert T(x) - x \right\vert^2 d \mu_0(x)
\end{equation}

It can be shown that such a mapping is the gradient of a convex function $u$, that is: $T(x) = \nabla_{x} u(x)$, which makes it irrotational (meaning $\nabla \times T = 0$~\cite{BuddWilliams}) if the underlying manifold is $\mathbb{R}^d$. More specifically, the regularity theory of Caffarelli (see the summary in~\cite{villani1}) shows us that if $f_0, f_1 \in C^{0, \alpha}(\mathbb{R}^d)$ for $0<\alpha <1$ and $0 < f_0, f_1 < \infty$, then we are guaranteed that the mapping is, in fact, a diffeomorphism.

In the more general manifold setting, we have that the mapping solving the Optimal Transport problem with squared geodesic cost is given by $T(x) = \text{exp}_{x} \left( \nabla u \right)$ (as long as the manifold is geodesically complete), see~\cite{mccann}, where $u$ is as a $c$-convex function, a technical definition that is common in the Optimal Transport literature, see~\cite{LoeperReg} for a discussion of this condition and more about the general case of manifolds. For the $n$-sphere, fairly general conditions on the source and target masses are given in~\cite{Loeper_OTonSphere}, which show that mass transports to a distance bounded strictly below the injectivity radius and gives us differentiability. Like the case of $\mathbb{R}^d$, we also need $0 <f_0, f_1 <\infty$. More detail on the regularity results of the Optimal Transport problem will be shown in~\ref{sec:optimaltransport}. The more general manifold case fails quite spectacularly in terms of differentiability, see the result in~\cite{figalli} for a simple example where differentiability fails to hold.

Optimal Information Transport provides another option for computing a mapping for moving mesh methods as a change of variables between two probability measures. We will provide much more detail on the Optimal Information Transport problem in Section~\ref{sec:optimalinformationtransport}, but we summarize the relevant details here. In this case, given a geodesic (in the space of probability measures endowed with the infinite-dimensional Fisher-Rao metric) between smooth source and target probability measures bounded away from zero, there exists a unique path $T_t: [0,1] \rightarrow \text{Diff}(M)$, where $\text{Diff}(M)$ denotes the space of diffeomorphisms over a manifold $M$. The key here is the manifold $M$ can be any compact, connected and complete manifold.

The result for Optimal Information Transport is geometrically more general than the corresponding result for Optimal Transport. In the case of the latter, (in the space of probability measures endowed with the Wasserstein metric) between smooth source and target probability measures bounded away from zero, there exists a unique path $T_t: [0,1] \rightarrow \text{Diff}$, where $\text{Diff}(M)$ in the cases $M = \mathbb{S}^d$ and $M = \mathbb{R}^d$. For some compact, connected and complete manifolds, again, as shown by the results in Figalli, et al~\cite{figalli}, this does not hold.

For the moving mesh problem, then, we expect both Optimal Transport and Optimal Information Transport to be versatile and useful in $\mathbb{R}^d$ and $\mathbb{S}^d$. Beyond these cases, Optimal Information Transport is preferable.

\section{Optimal Transport on Compact Manifolds}\label{sec:optimaltransport}
In this section, we consider the Monge problem of Optimal Transport on smooth compact, connected manifolds $M$ without boundary. Optimal Transport allows for a wide variety of cost functions to be chosen. In fact, the convergence framework in~\cite{HT_OTonSphere} allows for a range of important cost functions. However, in this manuscript we will focus on the squared geodesic cost, that is $c(x,y) = \frac{1}{2} d_{M}(x,y)$. Given source and target probability densities $\mu_0$ and $\mu_1$ supported on the manifold that have density functions $f_0$ and $f_1$, respectively, that is $d \mu_0(x) = f_0(x) \text{dvol}(x)$ and $d \mu_1(x) = f_1(y) \text{dvol}(y)$, find a mapping $T:M \rightarrow M$ such that $T_{\#} \mu_0 = \mu_1$, that is $\mu_1$ is the pushforward measure of $\mu_0$ via the mapping $T$ such that:

\begin{equation}
T = \argmin_{S} \int_{M} \frac{1}{2}d_{M} \left(x,S(x) \right) f_0(x)dx
\end{equation}

The Monge problem of Optimal Transport can be formulated as a partial differential equation given some additional assumptions. The squared geodesic cost function satisfies the following hypotheses~\ref{hyp:costfunction}, which are derived from the hypotheses given in~\cite{MTW}:

\begin{hypothesis}[Conditions on cost function]\label{hyp:costfunction}
The squared geodesic cost function $c(x,y) = \frac{1}{2} d_{M}(x,y)$ satisfies the following conditions:
\begin{enumerate}
\item[(a)] The function $\nabla_{M, x} c(x,y)$ is Lipschitz for all $x, y \in M$ and invertible in $y$
\item[(b)] For any $x, y \in M$, we have $\det D^2_{M, xy} c(x,y) \neq 0$
\end{enumerate}
\end{hypothesis}

As first shown in~\cite{mccann}, the optimal map corresponding to this cost function is determined from the condition
\vspace*{-4pt}
\bq\label{eq:mapConditionSphere}
\begin{cases}
T(x,p) = \text{exp}_{x}(p), & x\in M, p\in \Tf_x\\
\end{cases} \vspace*{-4pt}
\eq
where $\Tf_x$ denotes the tangent plane at $x$. The solution to the optimal transport problem is then given by the fully nonlinear partial differential equation
\bq\label{eq:PDE1}
F \left( x,\nabla_{M} u(x), D_{M}^2u(x) \right) = 0
\eq
where
\bq\label{eq:OTPDE}
F(x,p,M) \equiv -\det \left( M+A(x,p) \right) + H(x,p)
\eq
subject to the $c$-convexity condition, which requires
\bq\label{eq:cconvex}
D_{M}^2u(x) + A(x,\nabla_{M} u(x)) \geq 0.
\eq
Here
\begin{align*}
A(x,p) &= D_{M,xx}^2c \left( x,T(x,p) \right)\\
H(x,p) &= \abs{\det{D_{M,xy}^2c \left( x,T(x,p) \right)}}f_0(x)/f_1 \left( T(x,p) \right),
\end{align*}
and the PDE now describes a nonlinear relationship between the surface gradient and Hessian on the manifold. For the case of the sphere, we have the simple formula derived in~\cite{HT_OTonSphere2}

\begin{equation}
H(x,p) = \left\Vert p \right\Vert f_0(x) / \sin \left\Vert p \right\Vert f_1\left( T(x,p) \right)
\end{equation}

\subsection{Regularity}
We consider the optimal transport problem~\eqref{eq:PDE1} under the following hypothesis.

\begin{hypothesis}[Conditions on data (smooth)]\label{hyp:Smooth}
We require problem data to satisfy the following conditions:
\begin{enumerate}
\item[(a)] There exists some $m>0$ such that $f_1(x) \geq m$ for all $x\in \mathbb{S}^n$.
\item[(b)] The mass balance condition holds, $\int_{\Sf} f_0(x)\,dx = \int_{\Sf} f_1(y)\,dy$.
\item[(d)] The data satisfies the regularity requirements $f_0, f_1 \in C^{1,1}(\mathbb{S}^n)$.
\item[(f)] The geometry is the $n$-sphere $\mathbb{S}^n$
\end{enumerate}
\end{hypothesis}

\begin{comment}
The second set of hypotheses~\ref{hyp:Nonsmooth} will lead to \textit{a priori} unique Lipschitz Alexandrov viscosity solutions~\cite{LoeperReg}. However, these Lipschitz Alexandrov solutions are, in fact, Lipschitz viscosity solutions, details are in Appendix~\ref{app:viscosity}.

\begin{hypothesis}[Conditions on data (non-smooth)]\label{hyp:Nonsmooth}
We require problem data to satisfy the following conditions:
\begin{enumerate}
\item[(a)] The mass balance condition holds, $\int_{\Sf} f_1(x)\,dx = \int_{\Sf} f_2(y)\,dy$.
\item[(c)] The cost function $c(x,y)$ is the squared geodesic cost or $f(d(x,y))$ for some $f: [0, \pi) \rightarrow \mathbb{R}$ smooth and strictly increasing with $f'(0) = 0$ and satisfying that for some $C>0$ and some $\epsilon>0$ that $f''(t) \leq C$ for all $t \leq \epsilon$.
\item[(d)] There exists some $m>0$ such that $f_2(x) \geq m$ for all $x\in \mathbb{S}^n$.
\item[(e)] The source mass density $f_1$ is upper semi-continuous and the target mass density $f_2$ is continuous.
\item[(f)] The source and target masses are bounded $f_1, f_2 < \infty$.
\end{enumerate}
\end{hypothesis}
\end{comment}

We thus get the following regularity result, from Loeper~\cite{Loeper_OTonSphere}

\begin{theorem}[Regularity]\label{thm:regularity}
The optimal transport problem~\eqref{eq:PDE1} with data satisfying Hypothesis~\ref{hyp:Smooth} has a classical solution $u\in C^3(\Sf)$.
\end{theorem}

The solution to~\eqref{eq:PDE1} is unique only up to additive constants.  In this manuscript, sometimes we will fix a point $x_0 \in M$ and add the additional constraint:
\bq\label{eq:fixed}
u(x_0) = 0.
\eq
At other times, it may be more convenient to choose the mean-zero solution, that is to impose the constraint
\bq\label{eq:meanZero}
\int_{M} u = 0.
\eq

\subsection{Beyond the Sphere}
In some cases, manifolds $M$ with non-positive cost-sectional curvature at any point $x \in M$ can have positive measures $\mu_0, \mu_1 \in C^{\infty}(M)$, but $T$ is not even guaranteed to be continuous~\cite{LoeperReg}. We introduce the Ma, Trudinger, Wang tensor~\cite{MTW, LoeperReg}:

\begin{equation}
\mathfrak{G}_{c}(x_0, y_0) (\xi, \nu) = D^4_{p_{\nu}, p_{\nu}, x_{\xi}, x_{\xi}} \left[ (x,p) \mapsto -c \left(x, T_{x_{0}}(p)  \right) \right] \vert_{x_{0}, p_{0} = - \nabla_{x} x(x_0, y_0)}
\end{equation}
The cost-sectional curvature is negative at a point $(x, y)$ if there exist $\xi, \nu$ such that $\mathfrak{G}_{c}(x,y)(\xi, \nu) < C_0 \left\vert \xi \right\vert^2 \left\vert \nu \right\vert^2$. The idea is more perhaps transparent if one looks at the squared geodesic cost $c(x,y) = \frac{1}{2}d_{M}(x,y)$. In this simple case, the non-positivity of the cost-sectional curvature is equivalent to the negativity of the sectional curvature of $M$ at any point $x \in M$ due to the equality:

\begin{equation}
\frac{\mathfrak{G}_{c}(x,x)(\nu, \xi)}{\left\vert \xi \right\vert^2 \left\vert \nu \right\vert^2 - (\xi \cdot \nu)^2} = \frac{2}{3} \cdot \text{sectional curvature of $M$ at $x$ in the plane $(\xi, \nu)$}
\end{equation}

But even if the underlying manifold has strictly positive curvature everywhere, the Ma, Trudinger, Wang tensor $\mathfrak{G}_{c}(x,y)$ may be negative on the off-diagonal $(x,y) \ \text{s.t} \ x \neq y$. This was shown true even for some ellipsoids of revolution in the paper~\cite{figalli}. This suggests that computing Optimal Transport for such cases is more difficult and is not guaranteed to result in a diffeomorphic mapping $T$.

\section{Optimal Information Transport on Compact Manifolds}\label{sec:optimalinformationtransport}
The Monge problem of Optimal Transport hints at a deeper relation between the space of probability distributions and the mappings between probability measures. We have seen in section~\ref{sec:optimaltransport} that for the squared geodesic cost $c(x,y) = \frac{1}{2}d_{M}(x,y)$ there exists a unique mapping $T$ of the Monge problem which is then used in computing the total cost of transporting the probability measure $\mu_0$ to $\mu_1$. Furthermore, for the squared geodesic cost, this total cost defines a Riemannian distance between probability measures, see~\cite{villani1, villani2}, usually referred to as the Wasserstein distance. Furthermore, a Wasserstein interpolation between the source mass $\mu_0$ and the target mass $\mu_1$, denoted by $\mu(t)$ can be uniquely defined from a convex combination of the identity map and the Optimal Transport map from $\mu_0$ to $\mu_1$, that is, by defining the path of maps $T(t) = (1-t) \text{Id} + t T$ we can construct an interpolation of measures $\mu(t) = T(t)_{*} \mu_0$, which defines a geodesic in Wasserstein space.

A different distance between probability measures arising primarily in statistical applications is the Fisher-Rao distance. This distance is defined as the Riemannian distance arising from the Fisher-Rao metric on the space of probability measures. The Fisher-Rao metric, in turn, is the Hessian of the Kullback-Liebler divergence between the probability measures $\mu$ and $\nu$, that is the quantity

\begin{equation}
\text{KL}(\mu, \nu):= \int_{M} \log \left( \frac{d\mu}{d\nu} \right) d\mu
\end{equation}
which measures the relative entropy between one probability measure $\mu$ and another $\nu$. Here, the Fisher-Rao metric is the infinite-dimensional version, first studied in~\cite{friedrich}. Typically in statistics it is given in the finite-dimensional setting where the probability measures can be parametrized by a finite parameter space which is a subset of $k$-dimensional Euclidean space $\mathbb{R}^k$.

Given the space of probability measures then equipped with the Fisher-Rao metric, one may ask the question if there is a mapping $T$ which characterizes the Fisher-Rao distance in the same way that the Optimal Transport mapping $T$ characterizes the Wasserstein distance and interpolation. The short answer is yes, and one can derive it by analyzing the geometric connection between the space of diffeomorphisms $\text{Diff}(M)$ and the space of positive smooth probability measures (volume forms) $\text{Dens}(M)$ on the underlying manifold $M$. That is, one can define a projection map $\pi: \text{Diff}(M) \rightarrow \text{Dens}(M)$ providing a principle bundle structure, introduced by Moser~\cite{moser}. This structure is then analyzed to show that, sure enough, a unique Optimal Information Transport map $T$ can be associated with the Fisher-Rao distance between $\mu_0$ and $\mu_1$ and define a unique geodesic interpolation as well, and a simple system of PDE that can solved for $T$.

For more information about Optimal Information Transport, see the background and more detail presented in~\cite{modin, modin2}. Here we present a brief summary of the derivation of the Optimal Information Transport PDE presented in the manuscript~\cite{modin}. For all the discussion in this section, we will assume that the underlying manifold $M$ is compact and connected and is equipped with the standard volume form $\text{vol}$. We introduce the Fr\'{e}chet Manifold of smooth volume forms over $M$ with total volume $\text{vol}(M)$:

\begin{equation}
\text{Dens}(M) = \left\{ \mu \in \Omega^{n}(M): \int_{M} \mu = \text{vol}(M), \mu>0 \right\}
\end{equation}
where the notation $\Omega^n(M)$ denotes the $n$-forms over the manifold $M$, using topology induced by the Sobolev seminorms. The infinite-dimensional Fisher-Rao metric:

\begin{equation}\label{eq:fisherrao}
G^{F}_{\mu}(\alpha, \beta) = \frac{1}{4} \int_{M} \frac{d \alpha}{d \mu} \frac{d \beta}{d \mu} d\mu
\end{equation}
is weak in the Fr\'{e}chet topology, but the geodesics have explicit formulas! Furthermore, there is an explicit formula for the distance which is simply given by an integral:

\begin{equation}
d_{f}(\mu_0, \mu_1) = \sqrt{\text{vol}(M)} \arccos \left( \frac{1}{\text{vol}(M)} \int_{M} \sqrt{\frac{\mu_0}{\text{vol}} \frac{\mu_1}{\text{vol}}} \text{vol} \right)
\end{equation}

We introduce another manifold with rich structure. The set of diffeomorphisms on $M$ is denoted $\text{Diff}(M)$ with topology induced by the Sobolev seminorms.  We have assumed $M$ compact, which makes the set $\text{Diff}(M)$ is a Fr\'{e}chet Lie group under the composition of maps. The Lie algebra of $\text{Diff}(M)$ is given by the space $\mathfrak{X}(M)$ of smooth vector fields. We equip this manifold with the information metric on $\text{Diff}(M)$, defined by:

\begin{equation}\label{eq:complicatedmetric}
G^{I}_{\varphi} (U,V) = \int_{M} g(\Delta u, v) \text{vol} + \lambda \sum_{i=1}^{k} \int_{M} g(u, \xi_i) \text{vol}\int_{M} g(v,\xi_i) \text{vol}
\end{equation}
where $\lambda>0$, $u = U \circ \varphi^{-1}$, $v = V \circ \varphi^{-1}$, $g$ is the underlying metric of the manifold, $\Delta$ is the Laplace-de Rham operator on the space of vector fields, defined by $\Delta u = - \left( \delta d u^{\flat} + d \delta u^{\flat} \right)^{\sharp}$, where $\sharp$ and $\flat$ are the usual musical isomorphisms, $d : \Omega^{k}(M) \rightarrow \Omega^{k+1}(M)$ is the exterior differential and $\delta : \Omega^{k}(M) \rightarrow \Omega^{k-1}(M)$ is the codifferential, and $\{ \xi_i \}_i$ is an orthonormal basis of the harmonic vector fields on $M$, that is those vector fields $\xi$ for which $\Delta \xi = 0$. This metric~\eqref{eq:complicatedmetric} is also weak in the Fr\'{e}chet topology, but its geodesics are well-posed.

Now, the diffeomorphic maps define a group action on the space of densities, via the pullback $\varphi^{*} \mu = \nu$. The volume preserving maps (with respect to $\varphi$): define the isotropy group:

\begin{equation}
\text{Diff}_{\mu} (M) := \left\{ \varphi \in \text{Diff}(M); \varphi^{*} \mu = \mu \right\}
\end{equation}

This group action is transitive, that is, given $\mu$ and $\nu$, there exists a diffeomorphism $\varphi$ such that $\varphi^{*} \mu = \nu$. Therefore, the manifold of densities $\text{Dens}(M)$ is isomorphic to the space of right cosets $\text{Diff}(M) / \text{Diff}_{\mu} (M)$, and hence the density $\mu \in \text{Diff}(M)$ parametrizes a partition of the space of densities. Furthermore, by defining the projection map:

\begin{equation}
\pi_{\mu}: \text{Diff}(M) \ni \varphi \mapsto \varphi^{*} \mu \in \text{Dens}(M)
\end{equation}
we can define the following principal bundle structure:

\begin{center}
\begin{tikzcd}[column sep=small]
\text{Diff}_{\mu}(M) \arrow[r, hook] & \text{Diff}(M) \arrow[d, "\pi_{\mu}"] \\
 & \text{Dens}(M) 
\end{tikzcd}
\end{center}

The goal is now to give a Riemannian structure to the principal fiber bundle using~\eqref{eq:fisherrao} and~\eqref{eq:complicatedmetric}. The upshot is that although $\pi$ is a submersion it is actually also a \textit{Riemannian} submersion with respect to $G^F$ and $G^I$, i.e. the metric $G^I_{\varphi}$ descends to $G^F$. That is:

\begin{equation}
G^I_{\varphi}(U,V) = G^{F}_{\pi_{\mu}(\varphi)} \left( T_{\varphi} \pi_{\mu} \cdot U, T_{\varphi} \pi_{\mu} \cdot V \right)
\end{equation}
Descending metrics ``descend" horizontal geodesics, that is if $\varphi(t)$ is a horizontal geodesic curve in $\text{Diff}(M)$, then $\mu(t):= \pi(\varphi(t))$ is a geodesic curve in $\text{Dens}(M)$.

The problem of diffeomorphic density matching is that given a path of densities $\mu(t)$, we desire to find the horizontal path of diffeomorphisms $\varphi(t)$ which project onto $\mu(t)$, that are also of minimal length with respect to $G^I$. That is, solve the exact density matching problem, that is find a $\varphi(t)$ such that:

\begin{equation}\label{eq:exact}
\begin{cases}
\varphi(0) = \text{id}, \\
\varphi^{*} \mu_0 = \mu(t) \\
\text{minimizing} \ \int_{0}^{1} G^I_{\varphi(t)} \left( \dot{\varphi}(t), \dot{\varphi}(t) \right)dt
\end{cases}
\end{equation}

From now on, we assume that $\text{vol} = \mu_0$, i.e. $\mu_0$ is the canonical volume form. We take the equation $\varphi^{*}(t) \mu_0 = \varphi^{*}(t) \text{vol} = \mu(t)$ and differentiate with respect to $t$:

\begin{equation}
\dot{\mu}(t) = \partial_{t} \left( \varphi(t)^{*} \text{vol} \right) = \varphi^{*} \text{div}_{\text{vol}} v(t)
\end{equation}
where $v(t) = \dot{\varphi} \circ \varphi^{-1}$. This can be rewritten as:

\begin{equation}
\dot{\mu}(t) = \text{div} \left( v(t) \right) \circ \varphi(t) \mu(t)
\end{equation}

From here we perform the Hodge-Helmholz decomposition for the vector field $v$, by writing $v = \text{grad} f + w$. It turns out that the Hodge-Helmholtz composition is orthogonal with respect to the information metric $G^I$ so the length of the path $\varphi(t)$ is minimal for $w = 0$. Therefore, we must solve the following Poisson equation for the curl-free term $f$:

\begin{equation}\label{eq:diffmatching}
\begin{cases}
\Delta f(t) = \frac{\dot{\mu}(t)}{\mu(t)} \circ \varphi(t)^{-1} \\
\dot{\varphi}(t) = \text{grad}\left( f(t) \right) \circ \varphi(t), \quad \varphi(0) = \text{id}
\end{cases}
\end{equation}

As noted before, the geodesics for the Fisher-Rao metric have \textit{explicit forms}! Thus, we can insert explicit forms for $\mu(t)$ and $\dot{\mu(t)}$. Define $W: \text{Dens}(M) \rightarrow C^{\infty}(M)$ by $\mu \mapsto \sqrt{\frac{\mu}{\text{vol}}}$, then the geodesics are given by:

\begin{equation}
[0,1] \ni t \mapsto \left( \frac{\sin\left( (1-t) \theta \right)}{\sin \theta}f_0 + \frac{\sin(t \theta)}{\sin \theta}f_1 \right)^2\text{vol}
\end{equation}
where

\begin{equation}
\theta = \arccos \left( \frac{\left\langle f_0, f_1 \right\rangle_{L^2}}{\text{vol}(M)} \right)
\end{equation}
and $f_i = W(\mu_i)$. The Optimal Information Transport problem is then given by

\begin{equation}\label{eq:OIT}
\begin{cases}
\varphi(0) = \text{id}, \\
\varphi_{*} \mu_0 = \mu(t)
\end{cases}
\end{equation}
where $\varphi_{*} \mu_0$ denotes the pushforward of $\mu_0$. The solution of this problem~\eqref{eq:OIT}, the diffeomorphic mapping $T$, is then given by $T = \varphi^{-1}(1)$, where $\varphi(1)$ solves Equation~\eqref{eq:diffmatching}. We also have a solution for the inexact problem, where the object is to find $\varphi$ that minimizes:

\begin{equation}
\sigma d^2 (\text{id}, \varphi) + \bar{d}^2(\varphi^{*} \mu_0, \mu_1)^2, \quad \sigma>0
\end{equation}
where

\begin{equation}\label{eq:inexact}
d^2(\varphi, \psi) := \inf_{\gamma(0) = \varphi, \gamma(1) = \psi} \int_{0}^{1} G_{\gamma(t)} \left( \dot{\gamma}(t), \dot{\gamma}(t) \right) dt
\end{equation}
and $\bar{d}$ denotes the distance function on $\text{Dens}(M)$ with Riemannian metric $\bar{G}$. Thus, the first term in~\eqref{eq:inexact} is a regularity term and the second term penalizes $\varphi$ if the pushforward of $\mu_0$ is too far from $\mu_1$. Let $T = \varphi(1)^{-1}$ solve the exact problem~\eqref{eq:exact} then let $s = \frac{1}{1+\sigma}$. The solution to the inexact problem is then given by $T_{s} = \varphi(s)^{-1}$.

%\subsection{Time-Dependent Monitor Function in Moving-Mesh Methods}

%For adaptive mesh methods, the beauty of equation~\eqref{eq:OIT} is twofold. First, the geodesics $\mu(t)$ are explicit, as mentioned previously. Secondly, and most importantly for adaptive meshes, the form of $\mu(t)$ is \textit{not actually specified} and so we are free to choose it, provided that it is differentiable (or at least Lipschitz) in $t$. A convenient choice for moving mesh methods is to choose $\mu(t)$ to be equal to the monitor function $\mathcal{M}(t)$. Then solving the system~\eqref{eq:OIT} allows us to find the mapping $\varphi(t)$ that changes the initial grid at any time $s \in [0,1]$ to a target grid at any other time $t \in [0,1]$, such that $s\leq t$. This allows for a prescription $\varphi(t)$ that tells us how to do adaptive mesh methods in real-time!

\section{Numerical Methods for Optimal Transport and Optimal Information Transport}\label{sec:numerics}

We assume that we are on a compact, connected two-dimensional $C^3$ manifold $M$. The convergence framework and schemes we describe here are inspired by a numerical scheme recently designed by the authors for the $2$-sphere, which is equipped with a proof of convergence to the physically meaningful solution of the Optimal Transport problem, see~\cite{HT_OTonSphere,HT_OTonSphere2} for details and analysis. The adaptation of the convergence framework from the sphere to $M$ requires little modification provided that the solution to the elliptic PDE in question is merely \textit{a priori} know to be at least Lipschitz. We note that, of course, there exist other numerical schemes for computing the Poisson solver on the sphere or on compact, connected, complete manifolds $M$. We are presenting our provably convergent schemes to compare the results of applying Optimal Transport and Optimal Information Transport to problems on the sphere, the results of which will be shown in~\ref{sec:implementation}.

The discretizations presented here must be consistent and monotone. We also get compactness from Lipschitz stability. We will have in mind that we are using finite-difference schemes, but actually the convergence framework is more general than this. For the sphere, in the manuscript~\cite{HT_OTonSphere2} finite-difference schemes were, in fact, constructed which satisfied the hypotheses of the convergence framework. We have in mind that, similarly, finite-difference schemes can be constructed which satisfy the hypotheses for more general manifolds $M$, but defer the explicit construction of such schemes to later work. The properties of consistency, monotonicity, and Lipschitz stability are sufficient to guarantee convergence in nonsmooth cases and can be used to build explicit convergence rates, an avenue of research that is being pursued in a more general case (for the Optimal Transport PDE) in another manuscript by the author.

\subsection{Modification of PDE}

The Optimal Transport problem involves solving a fully nonlinear elliptic PDE. The Optimal Information Transport problem involves solving a Poisson equation. Due to the fact that we are solving on a compact manifold, the convergence framework we present here hinges first on the reformulation of each PDE to include some control on its Lipschitz constant. Given the underlying elliptic PDE that we wish to solve (such as a Monge-Amp\'{e}re-type PDE for the Optimal Transport problem or a Poisson PDE for the Optimal Information Transport problem) denoted by

\begin{equation}\label{eq:before}
F(x, Du(x), D^2 u(x)) = 0
\end{equation}
we instead discretize the following modified PDE:

\begin{equation}\label{eq:modified}
\max \{ F(x, \nabla u(x), D^2 u(x), \left\Vert \nabla u \right\Vert - R \} = 0
\end{equation}
where $R$ is chosen to be strictly greater than the \textit{a priori} Lipschitz bound. After discretization, the Lipschitz bound imparts sufficient stability on the solution to prove overall uniform convergence of the discrete solution $u^h$ to the actual solution $u$, provided that we have constructed a consistent, monotone scheme. The equivalence of the solution of~\eqref{eq:modified} to~\eqref{eq:before} for the Optimal Transport case was presented in the manuscript~\cite{HT_OTonSphere}. Here we present the argument for the Poisson equation.

First, we introduce the notion of a continuous viscosity solution, since the compactness arguments are very clean. The arguments for uniform convergence for elliptic PDE also are more naturally extended to weaker cases using the notion of viscosity solutions.

\begin{definition}[Continuous Viscosity Solutions]\label{viscosity}
An upper (lower) semicontinuous function $u: M \rightarrow \mathbb{R}$ is a viscosity sub-(super)- solution of the PDE~\eqref{eq:before} if for every $x_0 \in M$ and $\phi \in C^{\infty}(M)$ such that $u - \phi$ has a local maximum (minimum) at $x_0$ we have:

\begin{equation}\label{eq:viscineq}
F \left( x_0, u(x_0), \nabla u(x_0), D^2 u(x_0)\right) \leq (\geq) 0
\end{equation}
A continuous function $u: M \rightarrow \mathbb{R}$ is a viscosity solution of~\eqref{eq:before} if it is both a sub- and super-solution.
\end{definition}
The notion of test functions is slightly different in the case of the Optimal Transport PDE, as defined in~\cite{HT_OTonSphere}. Note also that, critically, the condition~\eqref{eq:viscineq} can hold vacuously. We modify the Poisson equation by defining the operator

\begin{equation}\label{eq:modifiedPoisson}
F(x, u(x), \nabla u(x), D^2 u(x)) := \max \left\{ -\Delta u + f(x),\left\Vert u(x) \right\Vert - R \right\} = 0
\end{equation}
The function $f(x)$ satisfies the compatibility condition: $\int_{M} f(x) dx = 0$. We now state the equivalence theorem that shows that adding a term to the modified Poisson equation~\eqref{eq:modifiedPoisson} does not change the solution. The proof of this theorem is given in Appendix~\ref{app:Poisson}.

\begin{thm}\label{thm:equivalencePoisson}
For any sequence $\epsilon_n \rightarrow 0$, $\epsilon_n>0$ and for $R>0$ large enough, the Lipschitz viscosity solution $w_{\epsilon_{n}}$
\begin{equation}\label{eq:modifiedPoissonepsilon}
\max \left\{-\Delta w_{\epsilon_{n}}(x) +f(x), \left\Vert w_{\epsilon_{n}}(x) \right\Vert - R \right\} + \epsilon_{n} w_{\epsilon_{n}}(x)= 0,
\end{equation}
where $f(x)$ satisfies $\int_{M} f(x)dx = 0$, converges to the mean-zero Lipschitz viscosity solution $u$ of the Poisson equation:

\begin{equation}\label{eq:Poisson}
- \Delta u(x) + f(x) = 0
\end{equation}
\end{thm}
This now allows us to instead discretize the PDE~\eqref{eq:modifiedPoisson} and likewise for the Optimal Transport problem.

\subsection{Grid, Computational Neighborhoods, \& Discrete Derivatives}

Now, we describe the grid and discretization on the manifold $M$. We start with a finite set of $N$ grid points $\G \subset M$ that lie on the manifold. We let $d_{M}(x,y)$ denote the usual geodesic distance between points $x,y$ on the manifold $M$. To the grid $\G$, we associate a number $h$ which serves as a parameter indicating the overall spacing of grid points. More precisely,
\begin{equation}\label{eq:h}
h = \sup\limits_{x\in M}\min\limits_{y\in\G^h} d_{M}(x,y) = \bO\left(N^{-1/2}\right).
\end{equation}
In particular, this guarantees that any ball of radius $h$ on the sphere will contain at least one discretization point. The grid must sufficiently resolve the manifold and satisfy some technical regularity requirements which are not difficult in practice to achieve, as detailed in~\cite{HT_OTonSphere}. These regularity requirements ensure that the grid is more or less isotropic with respect to direction. Now we seek an approximation of the form
\bq\label{eq:fd}
F^h(x,u;f_1,f_2) = 0, \quad x \in \G
\eq

At this point, we will begin to more explicitly state the construction of the finite-difference scheme in the case of the sphere $\mathbb{S}^2$. We defer the full generalization to compact surfaces $M$ to future work. We now consider a fixed grid point $x_i \in \G$ and a grid function $u:\G\to\R$. To perform the computation at $x_i$, we begin by projecting grid points close to $x_i$ onto the tangent plane at $x_i$.  That is, we consider the set of relevant neighboring discretization points
\bq\label{eq:neighbourCandidates}
\Zf(x_i) = \left\{z = \text{Proj}(x;x_i) \mid x \in \G, d_{M}(x,x_i) \leq \sqrt{h}\right\}.
\eq
In order to minimize the consistency error of the overall scheme, the choice of the radius of the neighborhood is $\sqrt{h}$ is optimal, see~\cite{FroeseMeshfreeEigs}, which will allow the directional derivatives to eventually resolve every direction in the computational neighborhood.

The projection onto local tangent planes is accomplished using geodesic normal coordinates, which are chosen to preserve the distance from $x_i$ (i.e. $d_{M}(x,x_i) = \|x_i -  \text{Proj}(x;x_i)\|$). This prevents any distortions that a non-zero curvature of a manifold necessarily introduces to second-order derivatives. This allows for monotone approximation schemes to be more easily constructed. Beyond the case of the sphere, it is necessary to merely have a consistent approximation for local geodesic coordinates, since explicit formulas for geodesic normal coordinates do not exist in general.

The Optimal Transport PDE requires computing gradients and eigenvalues of Hessian matrices. In order to compute the gradient, for example, we need to search various directions $\nu$ in the neighborhood. Therefore, we will consider the following finite set of possible directions,
\bq\label{eq:directions}
V = \left\{\left\{(\cos(jd\theta),\sin(jd\theta)),(-\sin(jd\theta),\cos(jd\theta))\right\} \mid j=1,\ldots,\frac{\pi}{2d\theta}\right\},
\eq
where the angular resolution $d\theta = \dfrac{\pi}{2\floor{\pi/(2\sqrt{h})}}$.

\begin{figure}
	\includegraphics[width=0.75\textwidth]{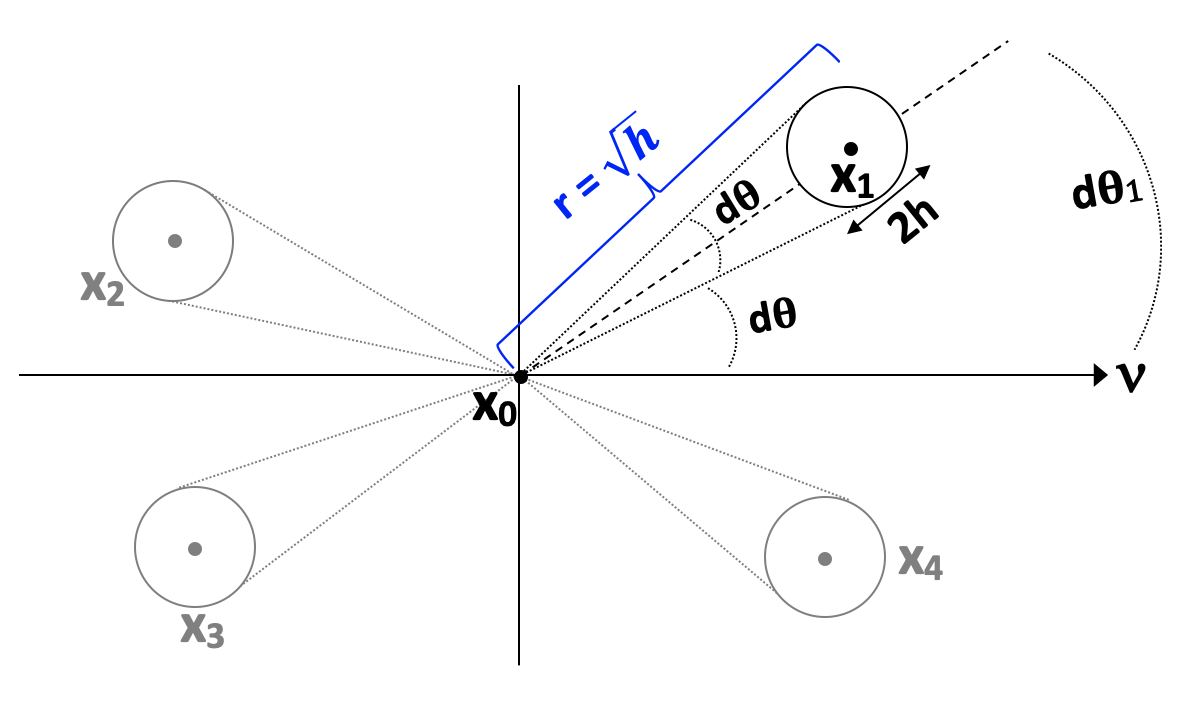}
	\caption{The computational points are chosen so as to ensure that they are within $2d\theta$ of the axis of computation, but also stay at least $d \theta$ away from the axis to ensure monotonicity of the scheme}\label{fig:neighborhood}
\end{figure}

In order to build approximations for derivatives in each direction $\nu\in V$, we will select four grid points $x_j\in\Zf(x_i)$, $j=1, \ldots, 4$, which will be used to construct the directional derivatives in this direction.  To accomplish this, we let $\nu^\perp$ be a unit vector orthogonal to $\nu$ and represent points in $x\in\Zf(x_i)$ using (rotated) polar coordinates $(r,\theta)$ centred at $x_i$ via
\[ x = x_i + r(\nu\cos\theta + \nu^\perp\sin\theta), x \in \Zf(x_i). \]
Then we select four points, each in a different quadrant ($Q_1, \ldots, Q_4$), that are well-aligned with the direction of $\nu$ and approximately symmetric about the axis $\nu$ via
\bq\label{eq:neighbours}
x_j \in \argmin\limits_{x\in\Zf(x_i)}\left\{\abs{\sin\theta} \mid \abs{\sin\theta} \geq d\theta, r \geq \sqrt{h}-2h, x\in Q_j\right\}
\eq
where $\cos\theta \geq 0$ for points in $Q_1$ or $Q_4$ and $\sin\theta \geq 0$ for points in $Q_1$ or $Q_2$, see Figure~\ref{fig:neighborhood}.

From here, we construct approximations of first directional derivatives (and second directional derivatives for the usual coordinate directions $(1,0)$ and $(0,1)$) of the form
\bq\label{eq:dnu}
\begin{split}
\Dt_{\nu}u(x_i) &= \sum\limits_{j=1}^4 b_j(u(x_j)-u(x_i)) \approx \frac{\partial u(x_i)}{\partial\nu}\\
\Dt_{\nu\nu}u(x_i) &= \sum\limits_{j=1}^4 a_j(u(x_j)-u(x_i)) \approx \frac{\partial^2u(x_i)}{\partial\nu^2}.
\end{split}
\eq

After performing a Taylor expansion of the sums in~\eqref{eq:dnu} we get explicit values for $a_j$ and $b_j$, which can be found in~\cite{HT_OTonSphere2, HT_OTonSphere3}. The discrete Laplacian is used in both our schemes for computing Optimal Transport and Optimal Information Transport. It can be built from approximations of second-order directional derivatives. Because the Laplacian is the trace of the Hessian, we take a geodesic normal coordinate system $(x^{j})_{j}$ and define:

\begin{equation}
\Delta^h u(x_i) = \mathcal{D}_{\nu_{x_{1}} \nu_{x_{1}}} u(x_i) + \mathcal{D}_{\nu_{x_{2}} \nu_{x_{2}}} u(x_i)
\end{equation}

In the case where the underlying PDE contains functions of the gradient, the goal is to add or subtract a Laplacian regularizing term in order to build monotonicity into the discretization, more details of which can be found in~\cite{HT_OTonSphere2}.

\subsection{Discretization of Optimal Transport and Optimal Information Transport PDE}
The full discretization for the discrete operator $F^h$ for the Optimal Transport problem on the sphere is:

\begin{multline}\label{eq:OTdisc}
F^h \left(x_i, u^h(x_i) \right) = \\
\min_{ \{ \nu_1, \nu_2 \} \in V} \prod_{j=1}^{2} \max \left\{ \mathcal{D}^h_{\nu_{j} \nu_{j}} u^h(x_i) + g^{-}_{1, \nu_{j}} \left( \nabla^h u^h(x_i) \right), 0 \right\} - f_1 (x_0) g^{+}_{2} \left( \nabla^h u^h(x_i) \right)
\end{multline}
where $g_1(p; \nu) = - \mathcal{D}^h_{\nu \nu} c(x_i,y) \vert_{y = T(x_i, p)}$ and $g_2(p; 
\nu) = \frac{\sin \left\Vert p \right\Vert}{\left\Vert p \right\Vert}$ and $g^{\pm}_{i, \nu} \left( \nabla^h u^h(x_i) \right) = g_{i} \left( \nabla^h u^h (x_i); \nu \right) \mp \epsilon_{g} \Delta^h u^h(x_i)$

And the full discretization for the discrete operator $F^h$ solving a Poisson equation on the sphere is:

\begin{equation}\label{eq:Poissondisc}
F^h \left(x_i, u^h(x_i) \right) = \Delta^h u^h(x_i) + \epsilon^h u^h(x_i) = f(x_i)
\end{equation}
where $\epsilon^h \rightarrow 0$ is chosen to be the consistency error of the operator $\Delta^h$. The discretization of the Lipschitz constraint is:

\begin{equation}
E^h \left( x_i, u^h(x_i) \right) := \max \left\vert \nabla^h u^h(x_i) \right\vert - R
\end{equation}
and we define:

\begin{equation}\label{eq:fulldisc}
G^h \left( x_i, u^h_n(x_i) \right) := \max \left\{ F^h \left( x_i, u^h_n(x_i) \right), E^h \left( x_i, u^h_n(x_i) \right) \right\}
\end{equation}

Then, the procedure for solving the discrete Optimal Transport PDE~\eqref{eq:OTdisc} is by using a parabolic scheme presented in Algorithm~\ref{alg:elliptic}:

\begin{algorithm}[h]
\caption{Computing the solution to elliptic PDE $F[u] = 0$}
\label{alg:elliptic}
\begin{algorithmic}[1]
\State Initialize $u^h_0$;
\State Fix $\epsilon>0$;
\While{$\left\vert G^h \left(x_i, u^h_n(x_i) \right) \right\vert > \epsilon$}
\State Compute $u^h_{n+1}(x_i) = u^h_n(x_i) + \Delta t G^h \left( x_i, u^h_n(x_i) \right)$
\EndWhile
\end{algorithmic}
\end{algorithm}

In the case of the Poisson equation when the Lipschitz constraint is inactive, we solve the linear system~\eqref{eq:Poissondisc} via standard linear algebra techniques.

\begin{remark}
The discretization~\eqref{eq:fulldisc} in principal involves verifying that the solution satisfies required Lipschitz bounds, then if necessary solving a second discrete problem to enforce the Lipschitz condition.  However, we do not see the verification step fail in practice, and hence never actually need to solve using the full operator $G^h$, but instead simply use the operator $F^h$.
\end{remark}

\subsection{Uniform Convergence Framework for the Poisson Equation for Compact Surfaces without Boundary}

The key properties in the discretization of the elliptic PDE without boundary for both the Optimal Transport and Optimal Information Transport PDE are consistency, monotonicity, and Lipschitz stability. Now, we show that with these properties we can derive a uniform convergence result for the Poisson equation on a compact 2D manifold without boundary. We use the auxiliary numerical scheme to prove boundedness:

\begin{multline}
\tilde{G}^h(v^h)(x_i) := \max \left\{ -\Delta^h v^h(x_i) - h^{\beta} v^h(x_i) + f(x_i), E^h \left( x_i, v^h(x_i) \right) \right\} + \\
h^{\alpha} v^h(x_i)= 0
\end{multline}

Let the formal consistency error of the scheme be $\mathcal{O}\left( h^{\alpha} \right)$.
Fix a point $x_0 \in M$. Choose the solution $u(x)$ such that $u(x_0) = 0$ and define $u^h$ such that $u^h(x_0) = 0$ by effecting the choice $u^h(x) = v^h(x) - v^h(x_0)$.

\begin{lemma}[Boundedness of Solution to Auxiliary Numerical Scheme]
The solution $v^h$ is bounded if $\beta \leq \alpha$.
\begin{proof}
The discretization $\tilde{G}^h$ is monotone and proper. Therefore, it has a discrete comparison principle, from the results in~\cite{ObermanSINUM}.
Then, we plug in $v^h$ and $u \pm C$ and apply the discrete comparison principle. The uniform bound follows from an identical argument for both the smooth and nonsmooth cases as shown in~\cite{HT_OTonSphere} for the Optimal Transport PDE.
\end{proof}
\end{lemma}
From the details of the proof in~\cite{HT_OTonSphere} we choose $\beta \leq \alpha$. It turns out that it is best to choose $\beta = \alpha$. We assume this choice from now on. We also assume that we can talk about $u^h(x)$, for any $x \in M$, assuming that we have used the interpolation of $u^h$ on the sphere detailed in~\cite{HT_OTonSphere}, with modifications, \textit{mutatis mutandis}, to apply the interpolation to the manifold $M$.

The uniform convergence critically depends on the solution of the modified PDE~\eqref{eq:modifiedPoisson} being the same as the solution of the Poisson equation. This was shown in~\ref{thm:equivalencePoisson}, while a similar result for the Optimal Transport PDE is in the paper~\cite{HT_OTonSphere}. Thus, we have:

\begin{thm}[Uniform Convergence in Nonsmooth Cases]
Assuming that we have a consistent and monotone discretization of $\Delta^h$ and the Lipschitz control on the discretization, then the solution $u^h$ converges uniformly to the Lipschitz viscosity solution $u$ of the Poisson equation on the manifold $M$ where $u(x_0) = 0$
\begin{proof}
Techniques are identical to those used in~\cite{HT_OTonSphere}. Uniform convergence can be proved for nonsmooth cases ($C^{0, 1}(M)$) by assuming that the modified Poisson equation~\eqref{eq:modifiedPoisson} has a unique solution.
\end{proof}
\end{thm}

\begin{remark}
The compactness arguments given in this section naturally do not lead to explicit convergence rates. However, forthcoming work will establish the guaranteed convergence rate $\left\vert u^h - u \right\vert = \mathcal{O} \left( h^{\alpha/3} \right)$ for \textit{a priori} $C^3$ solutions of the Poisson equation.
\end{remark}

\subsection{Algorithm for Optimal Information Transport}

In order to solve the Optimal Information Transport problem on $M$, we must solve the system of equations~\eqref{eq:diffmatching}. First, we must  perform some pre-computations involving the geometry and the mass densities. Then, we iterate an Euler scheme where each step involves one iteration of a Poisson solver. Here are the steps to compute the approximation to the diffeomorphic mapping $\varphi(1)$. The forward map at a time step $n$ will be denoted by $T_n$ and the corresponding inverse map at a time step $n$ will be denoted by $S_n$. We will outline the algorithm for the sphere, see Algorithm~\ref{alg:mapping}. More general manifolds require some fine-tuning, introduced in~\ref{sec:implementation}. The algorithm requires two important functions: $\text{Proj}$ a projection map consistent with the exponential map and $\text{Interp}$, a consistent interpolation map. For the case of the sphere, with an explicit formula for the exponential map, the function $\text{Proj}$ is simple to execute. However, the interpolation function is perhaps more challenging. A robust interpolation is that described in~\cite{HT_OTonSphere}, where the values are consistently interpolated over the Delaunay triangles.

\begin{algorithm}[h]
\caption{Computing the diffeomorphic mapping $\varphi = S(1)$}
\label{alg:mapping}
\begin{algorithmic}[1]
\State Initialize $T_0 = \text{id}$;
\State Initialize $S_0 = \text{id}$;
\State Fix $\Delta t \ll 1$;
\State Precompute $\theta$ via quadrature;
\While{$n \Delta t < 1$}
\State Compute the density function $\nu_{n} := \dot{\mu}_{n}/\mu_{n}$ using the explicit formulas;
\State Interpolate: $\text{Interp} \left\{ \nu_{n} \right\}$ onto the grid $\left\{ S(x_j) \right\}_{j}$;
\State Solve $\Delta^h f_n(x_i) = \nu_{n} \left( S_{n}(x_i) \right)$ with any monotone, consistent discretization of~\eqref{eq:modifiedPoisson};
\State Compute $\nabla^h f_n(x_i)$ for all $x_i$;
\State Interpolate: $\nabla^h f_n(x_i)$ onto the grid $\left\{ T(x_j) \right\}_{j}$;
\State Compute $T_{n+1}(x_i) = \text{Proj} \left\{ T_{n}(x_i) + \Delta t \nabla^h f_{n} \left( T_{n}(x_i) \right) \right\}$;
\State Compute $S_{n+1}(x_i) = S_{n} \left( \text{Proj} \left\{ x - \Delta t \nabla^h f_{n} (x_i) \right\} \right)$
\EndWhile
\end{algorithmic}
\end{algorithm}

\subsection{Adaptive Mesh Methods for compact, connected, and complete manifolds $M$}

\begin{figure}[htp]
	\subfigure[]{\includegraphics[width=0.49\textwidth]{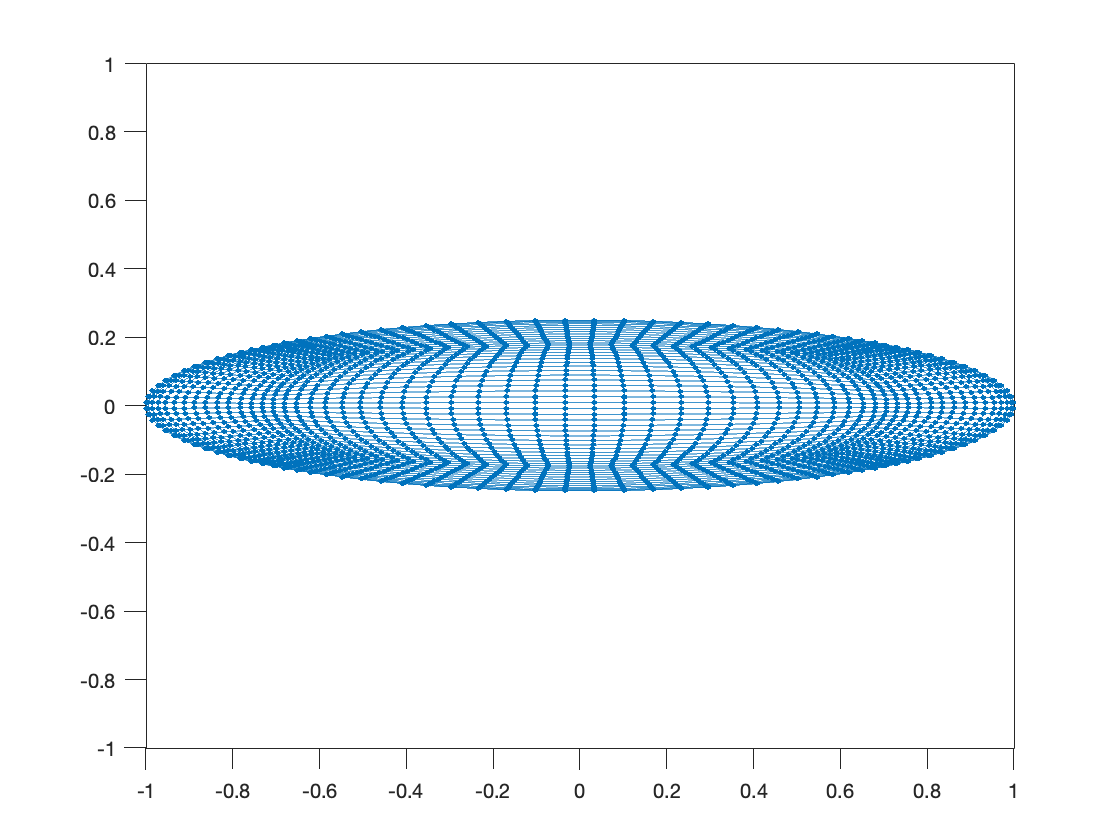}\label{fig:emeshbefore1}}
	\subfigure[]{\includegraphics[width=0.49\textwidth]{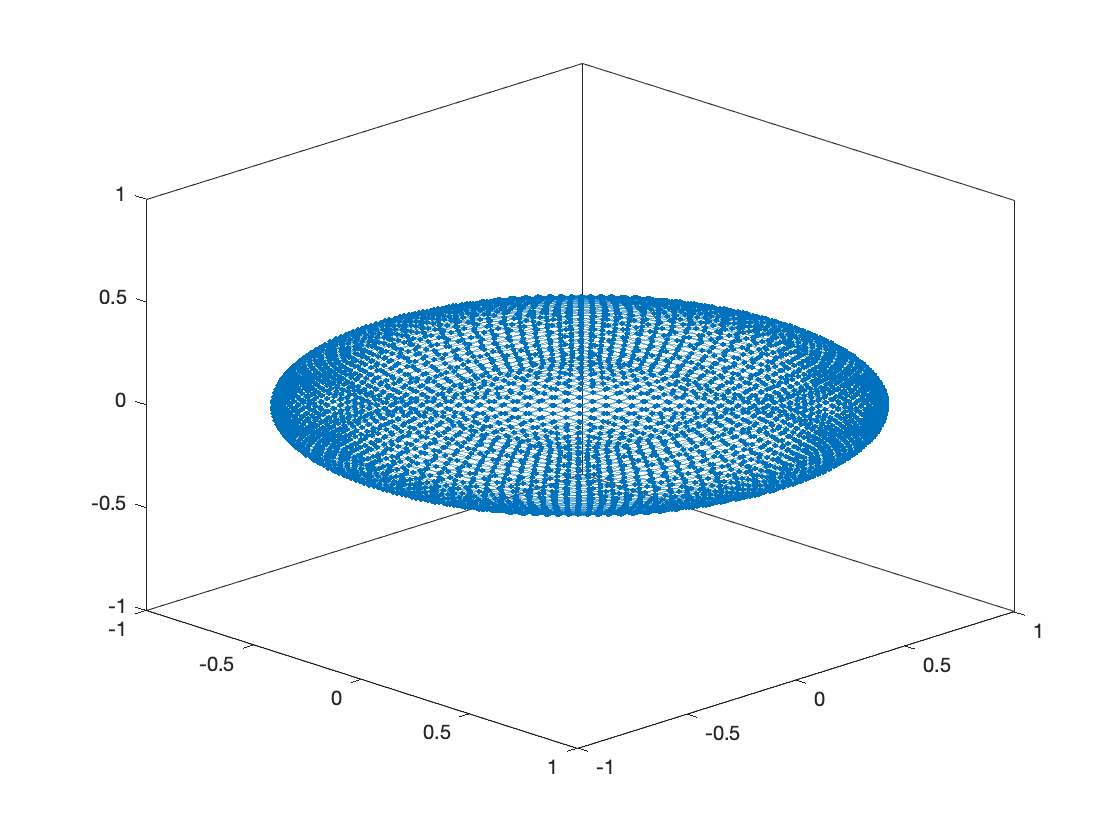}\label{fig:emeshbefore2}}
	\caption{A $N=5048$-point cube mesh for the ellipsoid from the side~\subref{fig:emeshbefore1} and diagonally above looking down~\subref{fig:emeshbefore2}}
	\label{fig:emeshbefore}
\end{figure}

Suppose, for example, that we have a mesh generated on an oblate sphere (ellipsoid of revolution with eccentricity $e = 0.25$), see~\ref{fig:emeshbefore}

For the Optimal Transport problem, the result from~\cite{figalli} asserts that there exist $C^{\infty}$ source and target mass density functions $f_0$ and $f_1$ such that the mapping $T$ is not guaranteed to be even continuous. Further difficulty in building a discretization comes again from the geometry regarding the computation of geodesics. However, given a local coordinate chart $(x^k)$ in a neighborhood of a point $x$, we can solve the second order ordinary differential equation for the geodesics:

\begin{equation}
\ddot{x}^k(t) + \dot{x}^i(t) \dot{x}^j (t) \Gamma^k_{ij} \left(x(t) \right) = 0
\end{equation}
where $\Gamma^k_{ij}$ are the Christoffel symbols. The Optimal Transport scheme we have presented, however, exploits an explicit representation formula for the mixed Hessian term, as derived in~\cite{HT_OTonSphere2}. This required an explicit formula for the exponential map as a surjective map from local tangent planes onto the sphere. This is no longer possible in more general geometries. Thus, a direct discretization of $D^2_{xy} c(x,y)$ would be necessary. We, similarly, do not have a representation formula for geodesic normal coordinates. This means that they must be approximated. This is possible in practice, however, with any method of finding geodesic normal coordinates that is consistent.

The only difficulties that exist with the method described for Optimal Information Transport are those concerning local geodesics. That means that this method is more easily generalizable than that of Optimal Transport. Furthermore, any provably convergent scheme for computing the solution of the Poisson equation on compact manifolds $M$ would lead to the same conclusions as we have presented here.

\section{Implementation}\label{sec:implementation}

Our objective is to perform adaptive mesh computations on the sphere. The original mesh will be assumed to be given by some off-the-shelf mesh generator which is simple to implement. Here we will be using a mesh defined on a cube which is then projected onto the sphere. The original mesh, therefore, \textit{does not have constant density}. However, in the moving mesh problem, we will be performing computations that produce diffeomorphic maps from a constant source density to a variable target density. This allows us to redistribute the mesh as desired.

\subsection{Successful Moving Mesh Method Implementation via Optimal Transport and Optimal Information Transport}

First, we demonstrate that the computation of mesh redistribution using Optimal Transport and Optimal Information Transport can be both successful in produce meshes which do not exhibit tangling. We select density functions with the goal of producing a transport map $T$ that will concentrate mesh points around the equator, see~\eqref{eq:meshgen}.

\begin{equation}\label{eq:meshgen}
\begin{cases}
f_1(\theta, \phi) = \frac{1}{4\pi} \\
f_2(\theta, \phi) = \left(1 - \text{exp} \left( -1/ 30 \left( \text{arccos}(z) - \pi/2 \right)^2 \right) \right)/3.53552
\end{cases}
\end{equation}

\begin{figure}[htp]
	\subfigure[]{\includegraphics[width=0.49\textwidth]{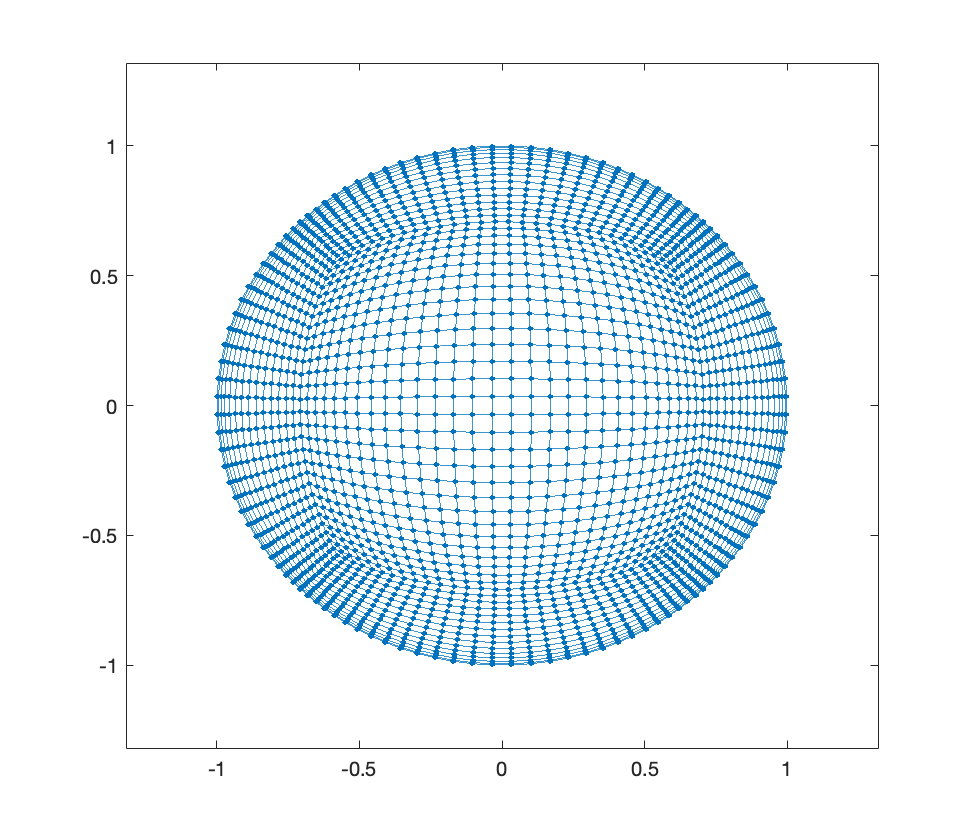}\label{fig:meshbefore1}}
	\subfigure[]{\includegraphics[width=0.49\textwidth]{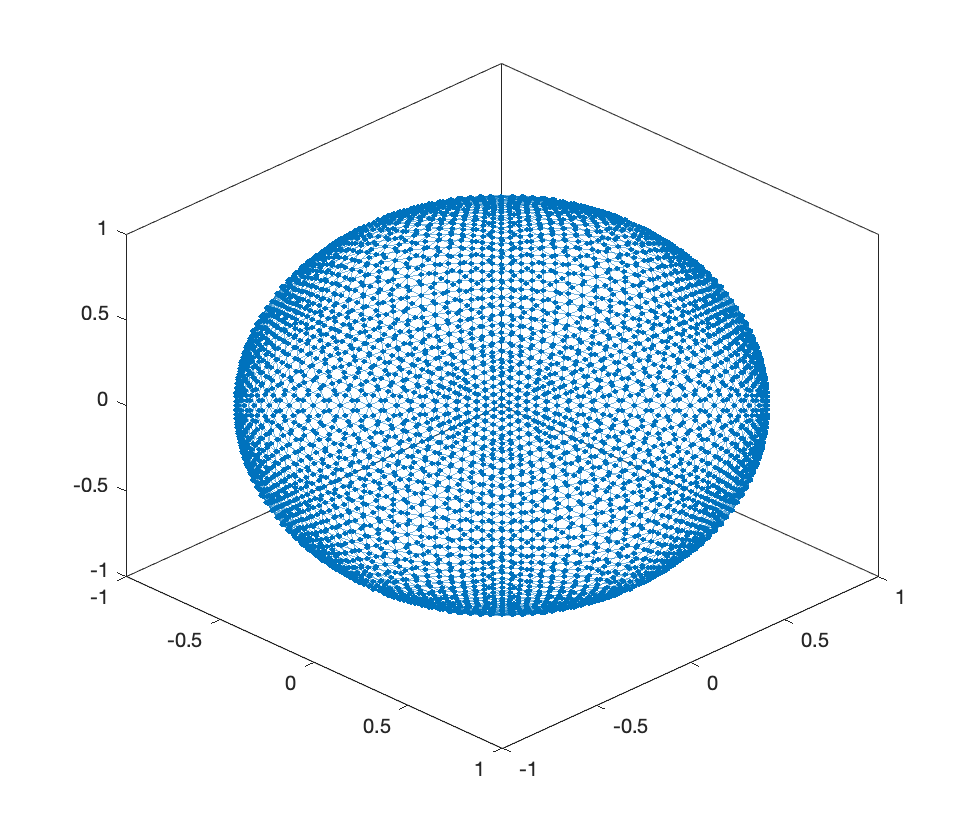}\label{fig:meshbefore2}}
	\caption{A $N=5048$-point cube mesh from the top~\subref{fig:meshbefore1} and diagonally above looking down~\subref{fig:meshbefore2}}
	\label{fig:meshbefore}
\end{figure}

The original mesh is generated from the cube mesh that is then projected onto the sphere, see~\ref{fig:meshbefore}

\begin{figure}[htp]
	\includegraphics[width=0.8\textwidth]{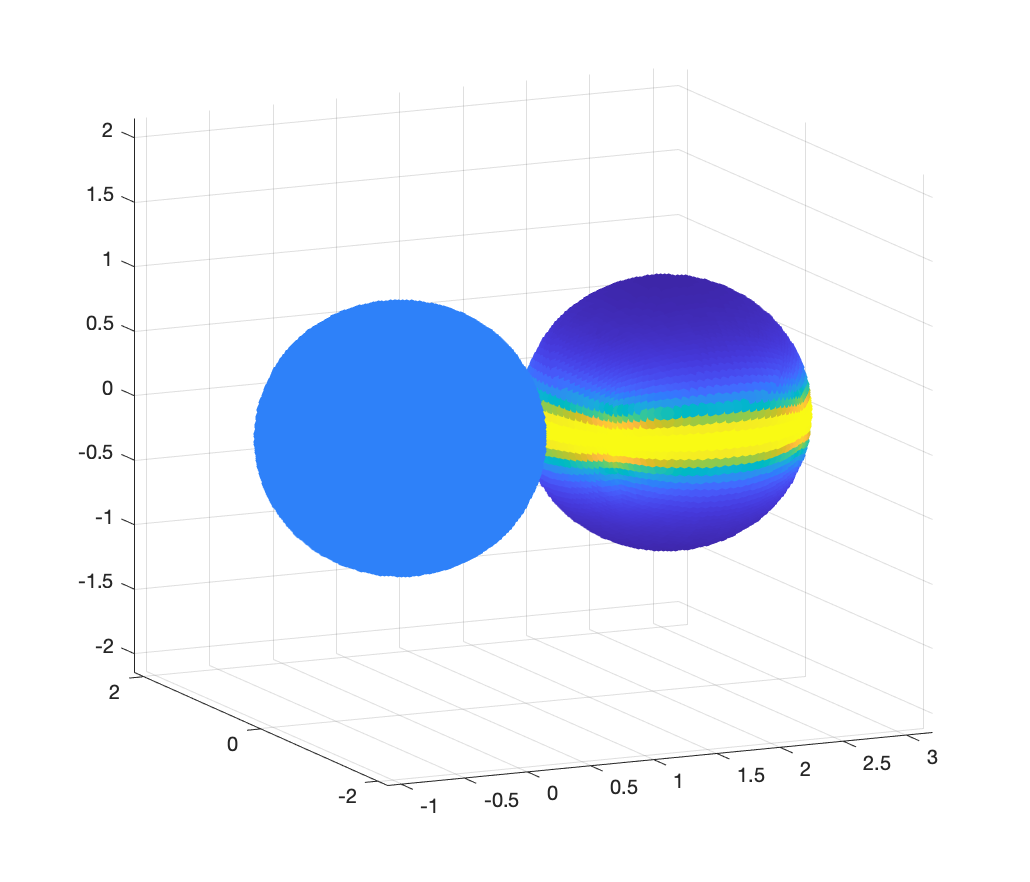}
	\caption{Density in the original mesh (left) and the modified density (right)}\label{fig:equator}
\end{figure}

We then require that the density concentrates about the equator according to~\eqref{eq:meshgen}, see the visual representation in the figure~\ref{fig:equator}

\begin{figure}[htp]
	\subfigure[]{\includegraphics[width=0.8\textwidth]{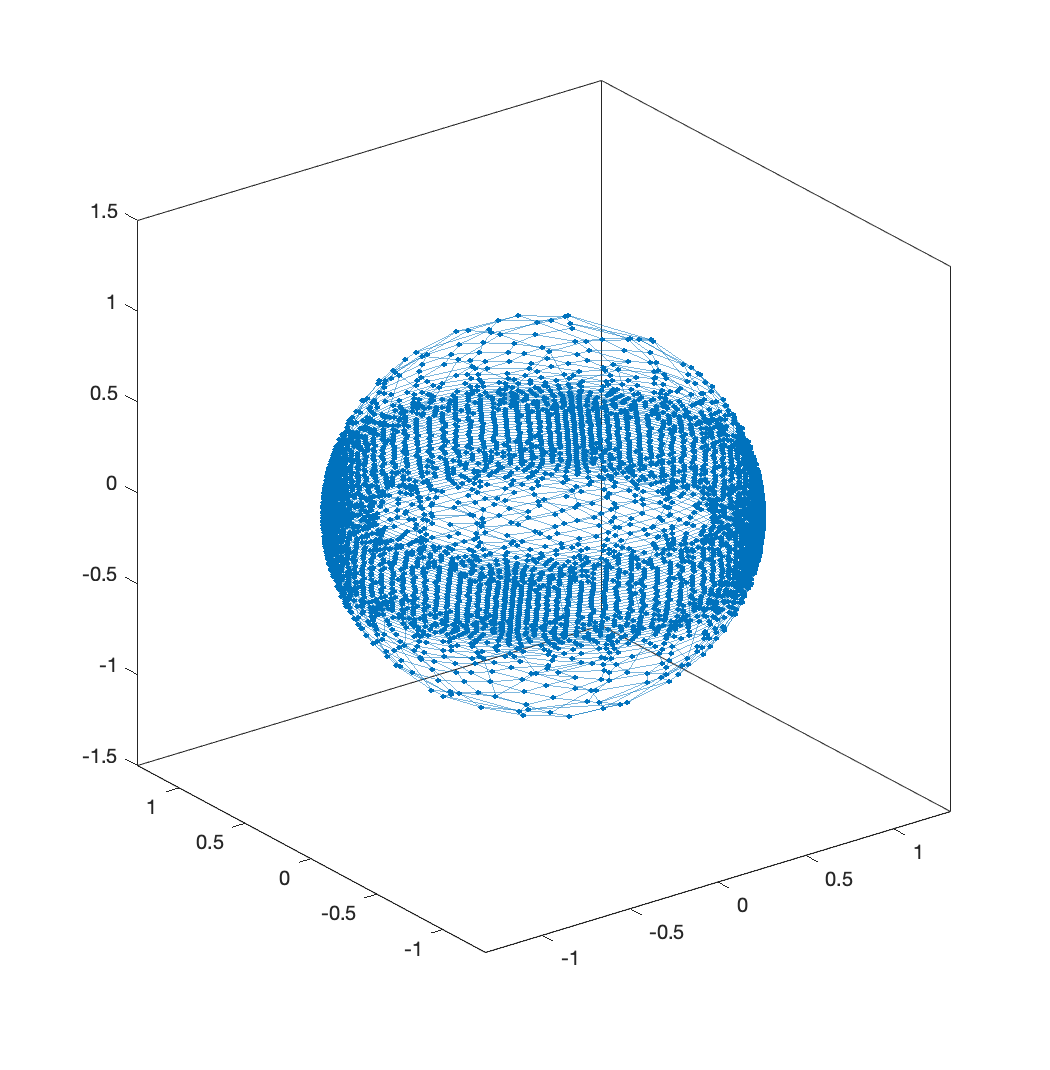}\label{fig:MAdiagonal}}
	\subfigure[]{\includegraphics[width=0.49\textwidth]{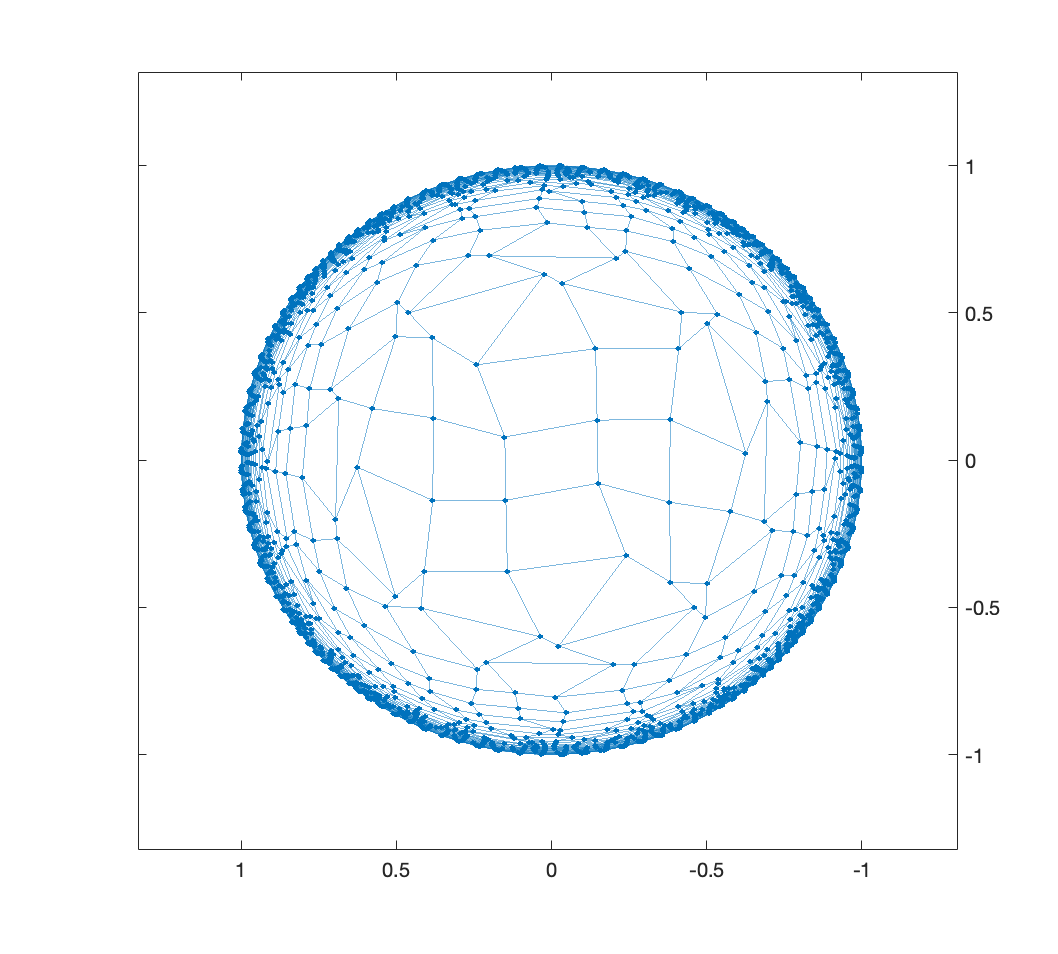}\label{fig:MAtop}}
	\subfigure[]{\includegraphics[width=0.49\textwidth]{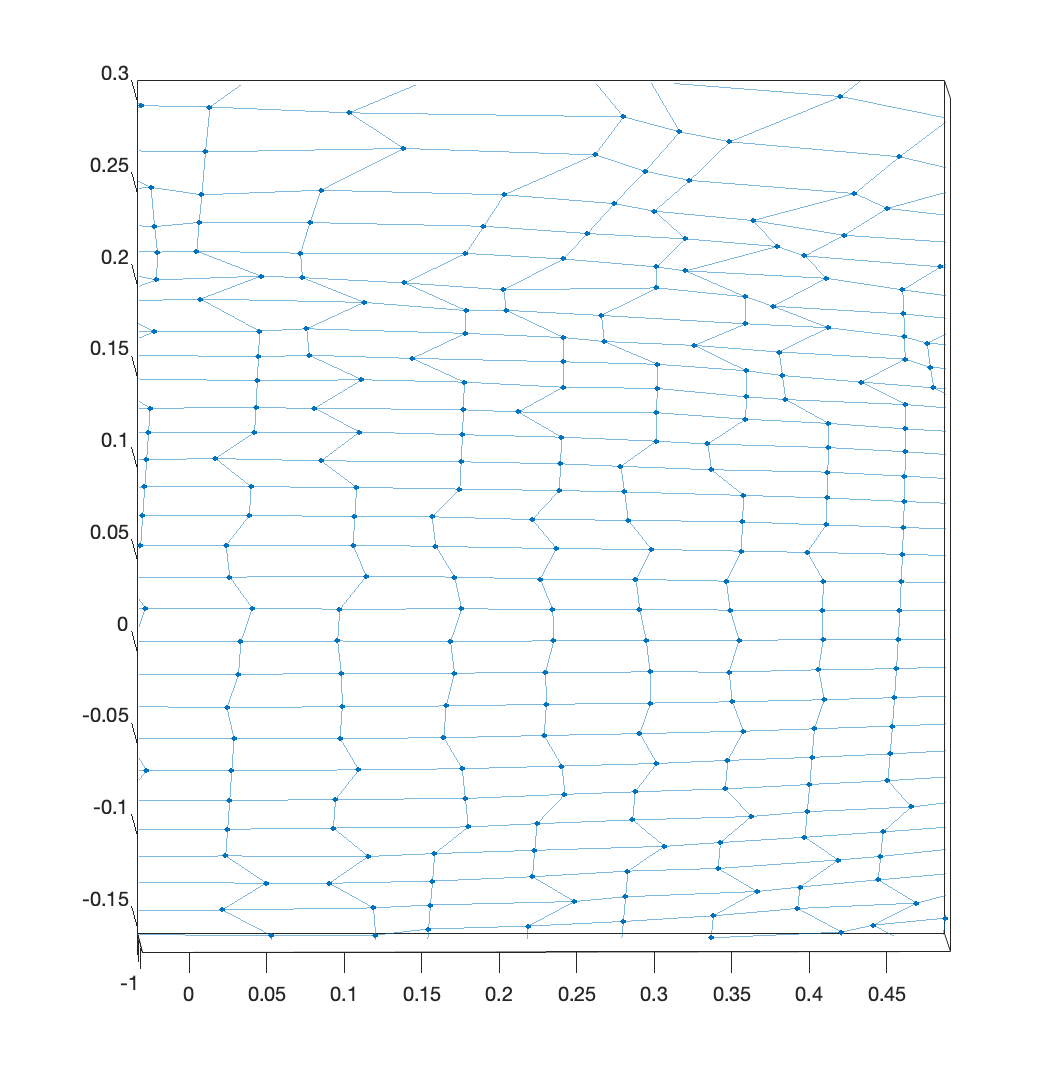}\label{fig:MAdetail}}
	\caption{~\subref{fig:MAdiagonal} 3D view of the mesh obtained via Optimal Transport map,~\subref{fig:MAtop} top view, and a detailed view of the equator~\subref{fig:MAdetail} showing that the grid lines do not tangle}
	\label{fig:MA}
\end{figure}

The resulting mesh restructuring computed using Optimal Transport is pictured in Figure~\ref{fig:MA}.

Now we perform the same computation using Optimal Information Transport, see~\ref{fig:movingmesh2}:

\begin{figure}[htp]
	\subfigure[]{\includegraphics[width=0.8\textwidth]{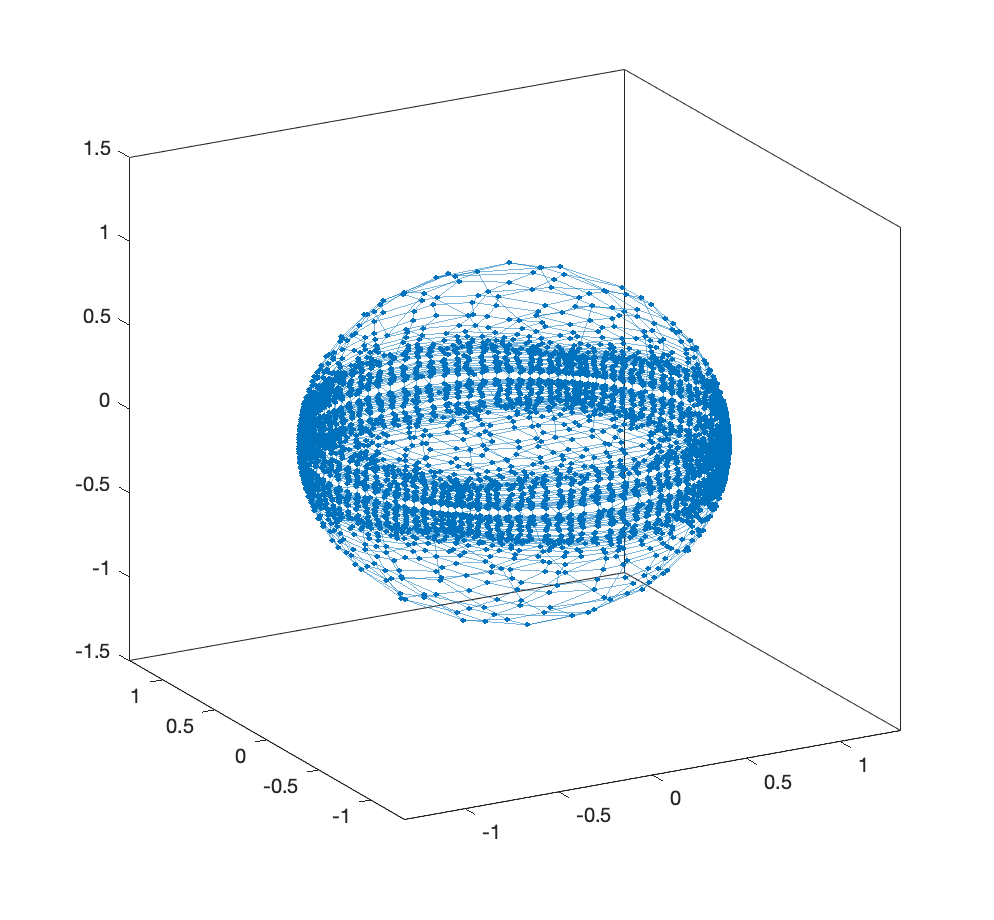}\label{fig:OITdiagonal}}
	\subfigure[]{\includegraphics[width=0.49\textwidth]{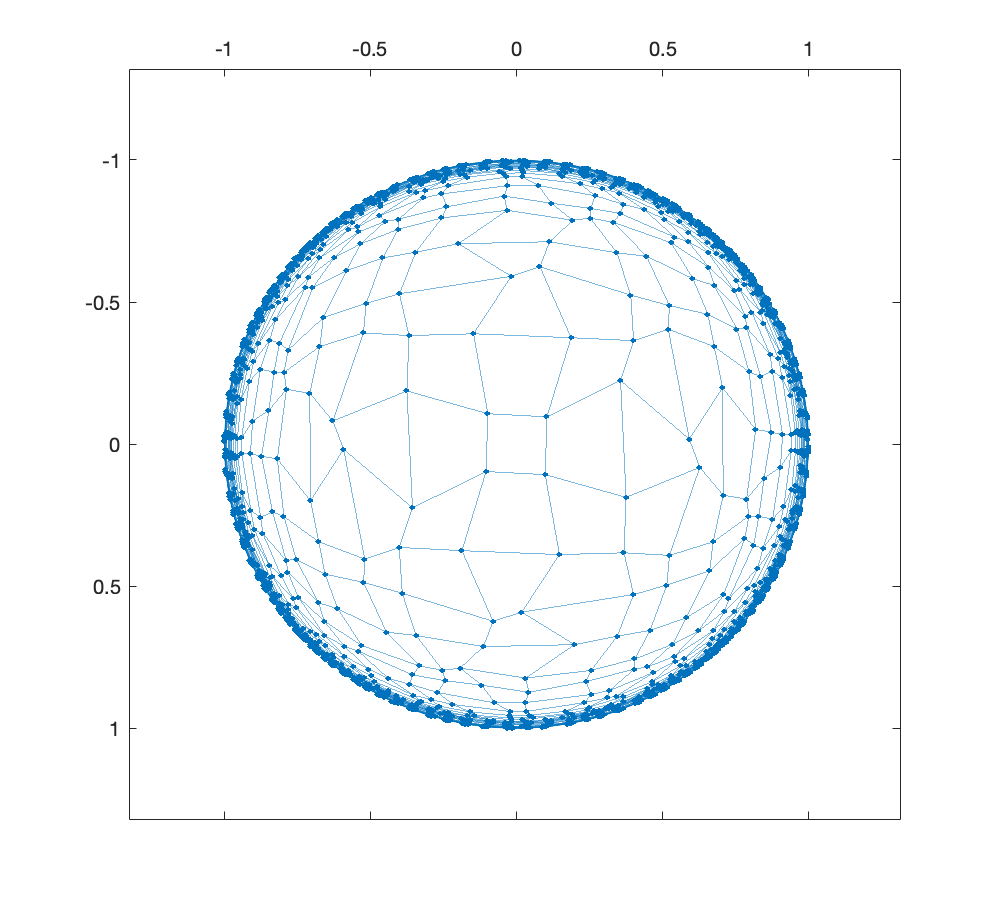}\label{fig:OITtop}}
	\subfigure[]{\includegraphics[width=0.49\textwidth]{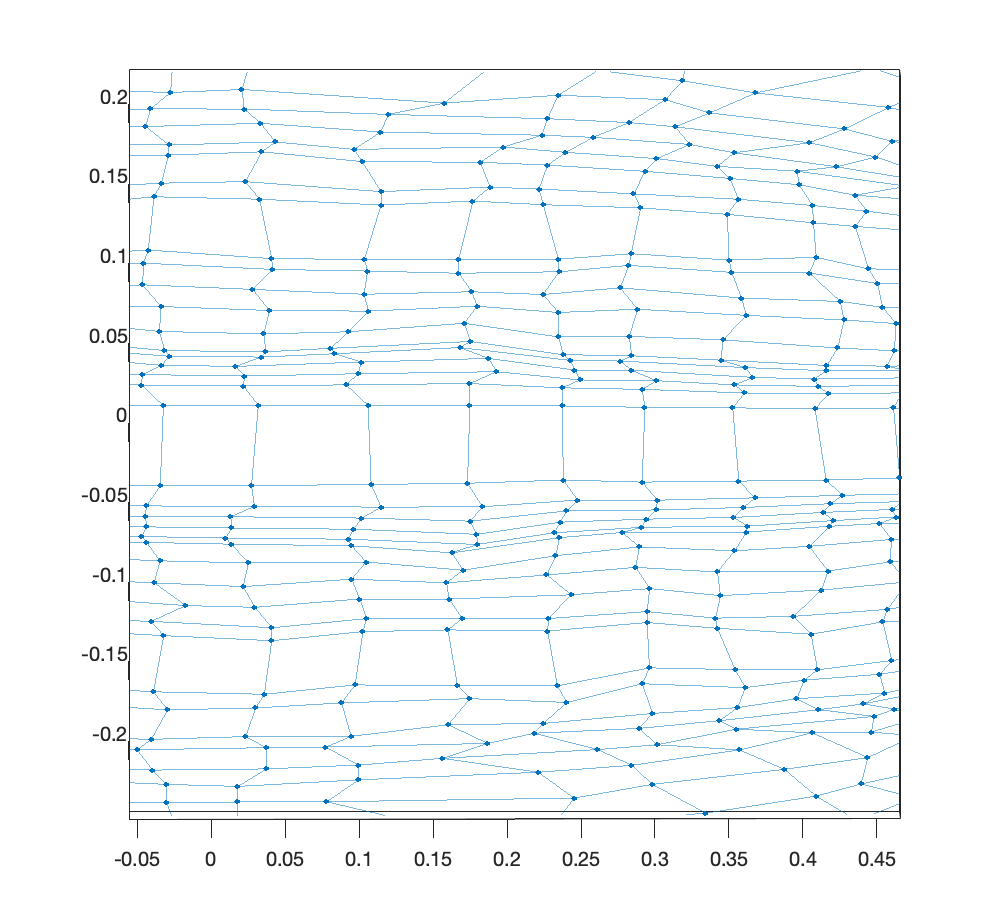}\label{fig:OITdetail}}
	\caption{~\subref{fig:OITdiagonal} 3D view of the mesh obtained via Optimal Information Transport,~\subref{fig:OITtop} top view, and a detailed view of the equator~\subref{fig:OITdetail} showing that the grid lines do not tangle}
	\label{fig:movingmesh2}
\end{figure}

\subsection{Optimal Information Transport is More Successful in Singular Examples}

\begin{figure}[htp]
	\includegraphics[width=0.8\textwidth]{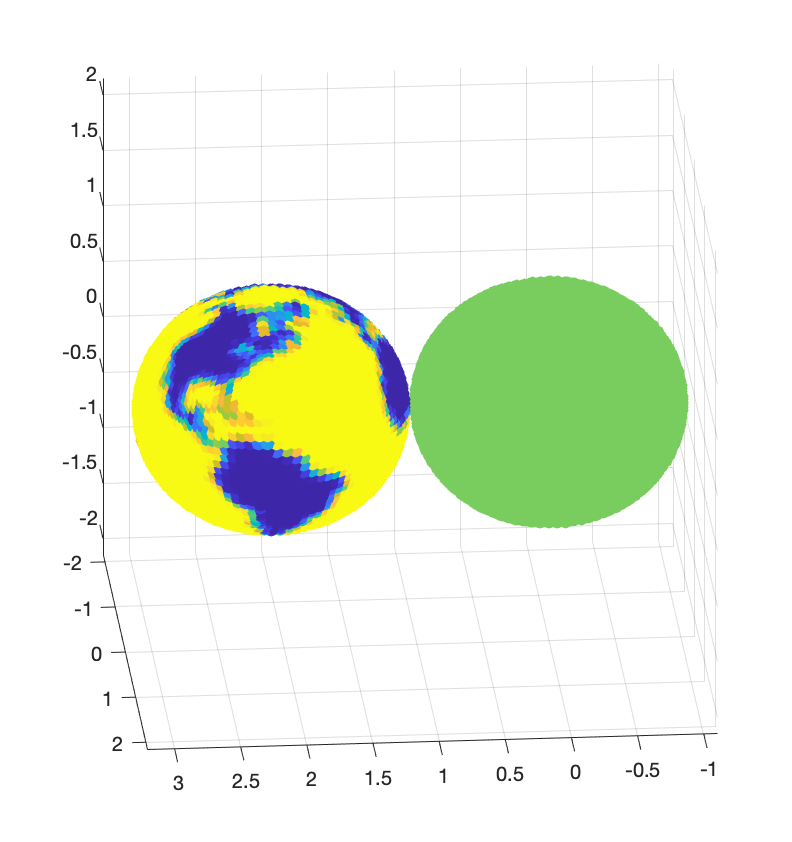}
	\caption{Constant density in the source density (right) and the target globe density (left)}\label{fig:fandgglobe}
\end{figure}

Another example serves to demonstrate the versatility and desirability of using Optimal Information Transport for moving mesh problems. We take a constant source density $f_1 = 1/(4 \pi)$ and map it to a complicated, discontinuous map of the world, where a higher density of grid points is required over the oceans and lower density over the continents, see Figure~\ref{fig:fandgglobe}.

\begin{figure}[htp]
	\subfigure[]{\includegraphics[width=0.49\textwidth]{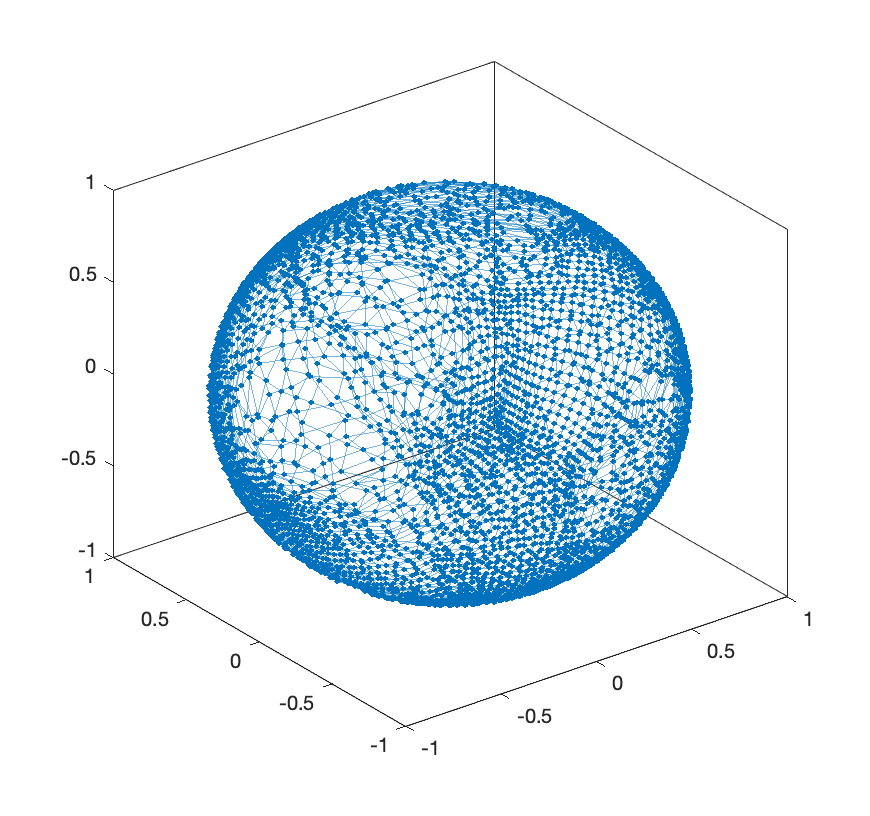}\label{fig:globediagonal}}
	\subfigure[]{\includegraphics[width=0.49\textwidth]{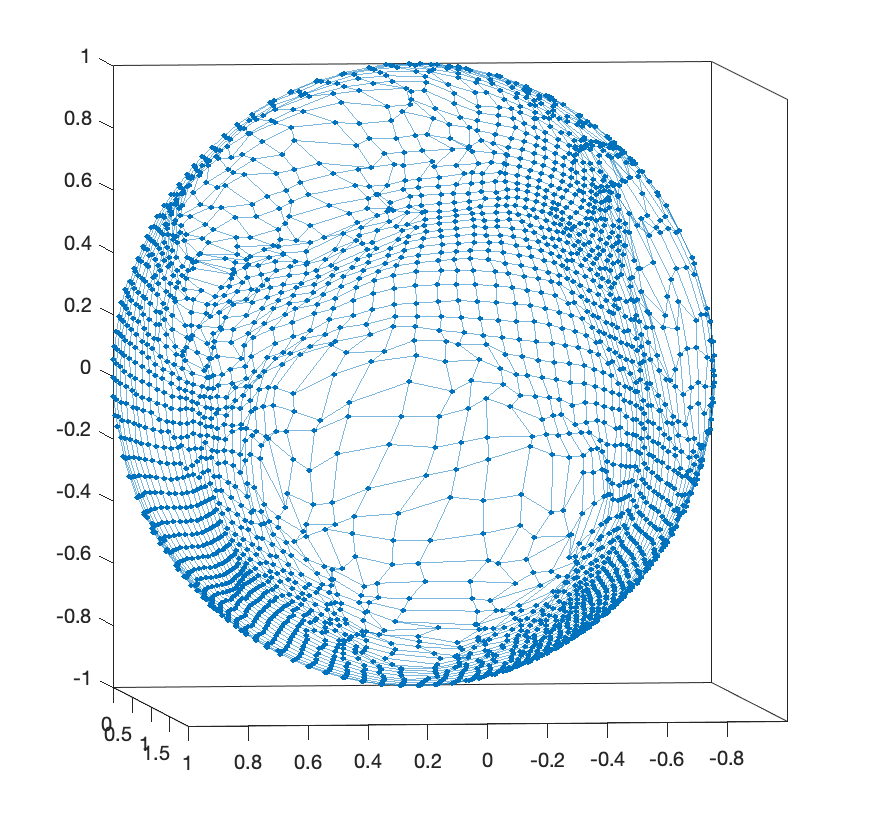}\label{fig:globeNA}}
	\caption{Entire spherical mesh transformed by Optimal Information Transport~\subref{fig:globediagonal} in 3D view. The view is of Africa, Europe and Asia, but the other side of the sphere mesh is also visible. Shown on the right is detail on North America with other side of sphere hidden from view in order to show the mesh more clearly~\subref{fig:globeNA}. The outlines of North and South America can be discerned and tangling is minimized.}
	\label{fig:globe}
\end{figure}

As before, we use the cube mesh~\ref{fig:meshbefore} and transform it via Optimal Information Transport. The result is shown in Figure~\ref{fig:globe}.

\begin{figure}[htp]
	\subfigure[]{\includegraphics[width=0.49\textwidth]{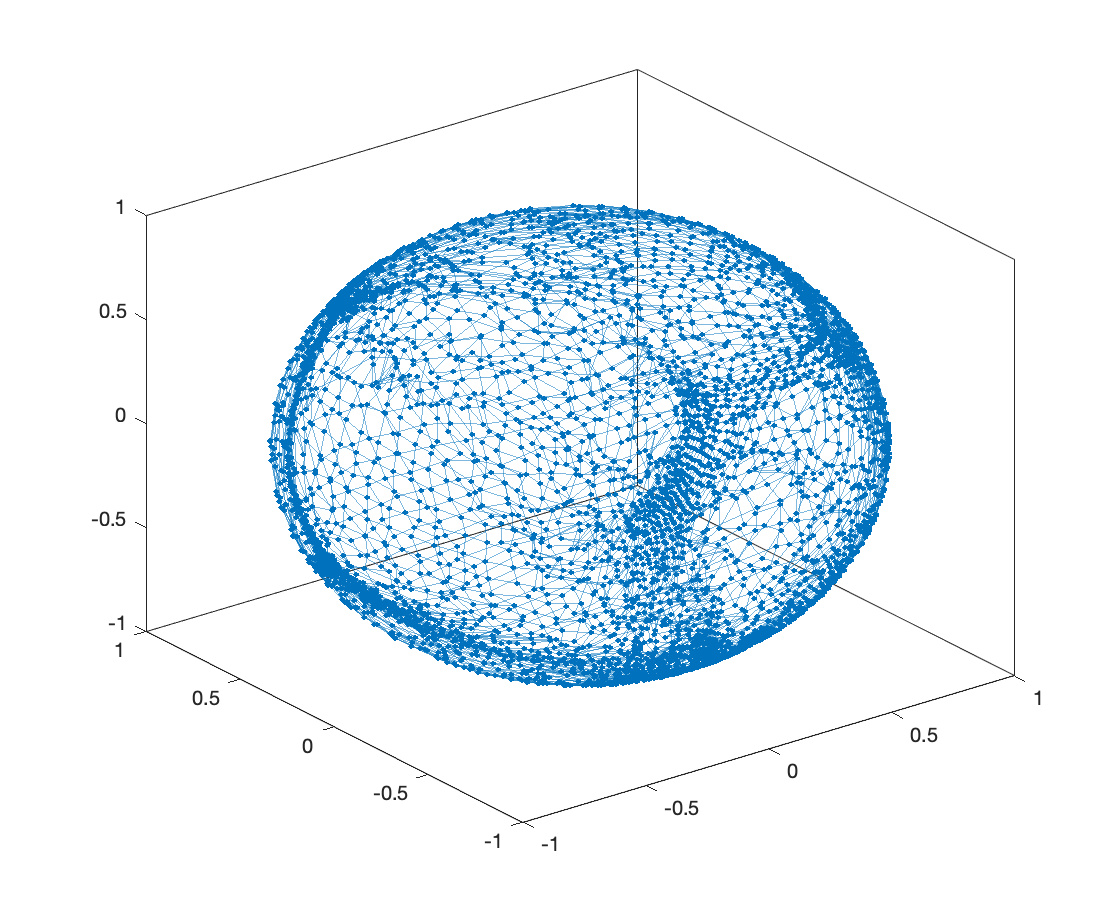}\label{fig:MAglobediagonal}}
	\subfigure[]{\includegraphics[width=0.49\textwidth]{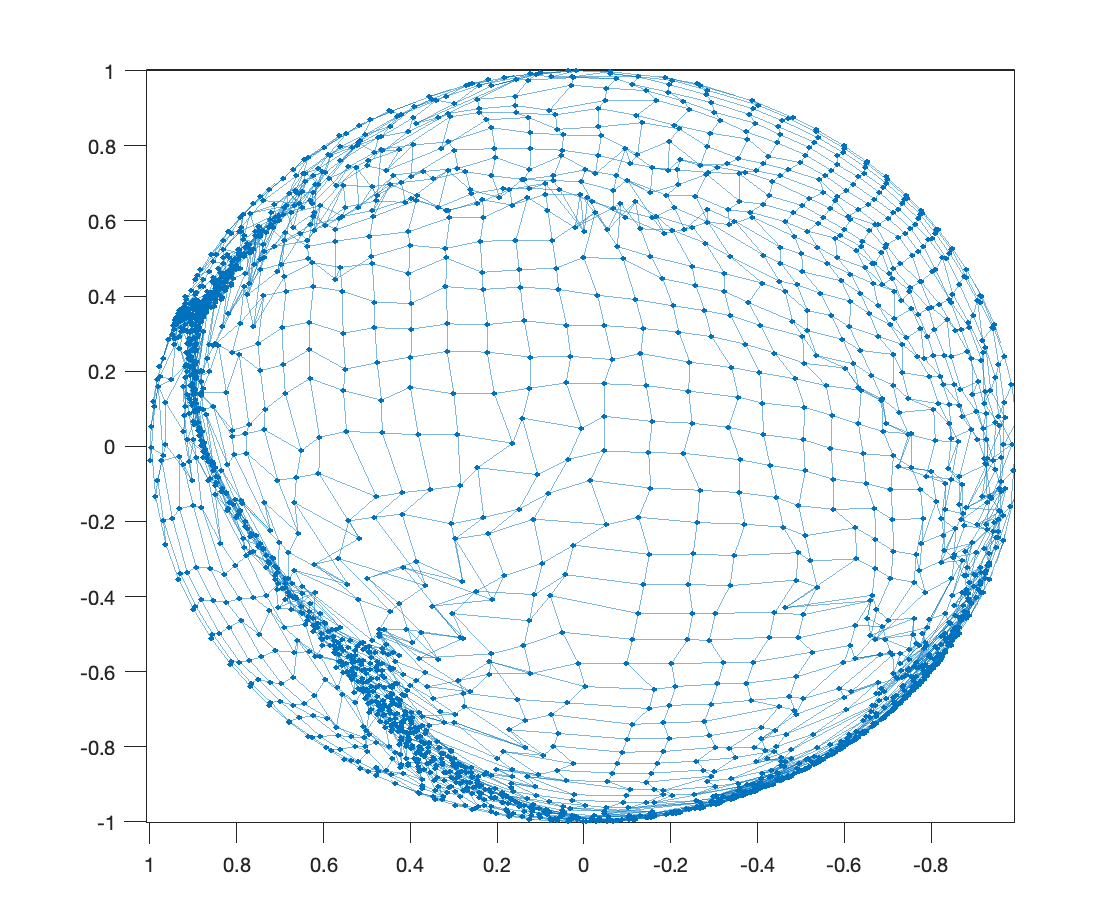}\label{fig:MAglobe}}
	\caption{Entire spherical mesh transformed by Optimal Information Transport~\subref{fig:globediagonal} in 3D view. Shown on the right is detail on Africa, the Middle East, and Asia with other side of sphere hidden from view in order to show the mesh more clearly~\subref{fig:MAglobe}. One can observe some tangling here.}
	\label{fig:globeMA}
\end{figure}

Some difficulties are encountered when attempting to perform the same numerical computation via Optimal Transport. First, the necessity of using an interpolation function for the target mass density $g$ onto the mesh defined by $T(x)$ leads to instabilities in the scheme due to the very nonlinear nature of the equation we are trying to solve. Compensating for this, by using an inverted image for the source mass $f$ and a constant for the target mass $g$ will run, but exhibits significant tangling for very complicated images like the world map, see Figure~\ref{fig:globeMA}

One potential solution to this issue exhibited by the Optimal Transport problem will be explored in further work where the computation of the mapping $T$ is slightly modified to gain convergence guarantees and greater stability.

\section{Conclusion}\label{sec:conclusion}
We have provided the first comparison of Optimal Transport and Optimal Information Transport on the sphere for the diffeomorphic density matching problem with applications to adaptive mesh methods by using convergent finite-difference schemes. We also introduced the issues concerning their generalization to more general closed, compact surfaces $M$. In this more general case, there appear to be more merits to using Optimal Information Transport. These merits include speed of implementation, generalizability in computational terms but also in terms of \textit{a priori} regularity results. Nevertheless, it could be that a particular application requires the computation of Optimal Transport with the squared geodesic cost on the sphere and more general surfaces.

As an additional benefit, the numerical techniques presented here are more computationally efficient than direct computation via discretization of the probability measures via Dirac measures. Many of the numerical techniques presented here apply more generally to elliptic PDE on compact closed surfaces, potentially solving a wide variety of problems.

\bibliographystyle{plain}
\bibliography{OTonSphere5}

{\appendix

\section{Modified Poisson}\label{app:Poisson}

Recall the definition of the viscosity solution~\ref{viscosity}.

\begin{thm}\label{thm:equivalencePoisson}
For any sequence $\epsilon_n \rightarrow 0$, $\epsilon_n>0$ and for $R>0$ large enough, the viscosity solution $w_{\epsilon_{n}}$ of the PDE
\begin{equation}\label{eq:modifiedPoissonepsilon}
\max \left\{-\Delta w_{\epsilon_{n}}(x) +f(x), \left\Vert \nabla w_{\epsilon_{n}}(x) \right\Vert - R \right\} + \epsilon_{n} w_{\epsilon_{n}}(x)= 0,
\end{equation}
and for continuous $f(x)$ satisfying $\int_{M} f(x) dx = 0$, converges uniformly to the unique mean-zero Lipschitz solution $u$ of the PDE $-\Delta u(x) + f(x) = 0$.
\begin{proof}
The equation~\eqref{eq:modifiedPoissonepsilon} is proper, and therefore even though it is degenerate elliptic, there exists a comparison principle, provided that $f(x)$ is continuous, by~\cite{CIL}. Thus, by Perron's method for continuous viscosity solutions, there exists a unique solution $w_{\epsilon}$ to this equation for any $\epsilon>0$.

By the same reasoning, the viscosity solution $u_{\epsilon}$ of the PDE $- \Delta u_{\epsilon} + f(x) + \epsilon u_{\epsilon} = 0$ is unique as well and is Lipschitz for all $\epsilon > 0$ provided that it remains uniformly bounded in $\epsilon$. It is mean zero, that is $\int_{M} u_{\epsilon} = 0$. The distributional solution $\tilde{u}_{\epsilon}$ of $- \Delta \tilde{u}_{\epsilon} + f(x) + \epsilon \tilde{u}_{\epsilon} = 0$ is the same as the viscosity solution $u_{\epsilon} = \tilde{u}_{\epsilon}$, see~\cite{ishii}. Thus, in the sense of distributions, for any $v \in C^2$ and use the integration by parts formula on the closed compact manifold $M$:
\begin{equation}
\int_{M} v \Delta u_{\epsilon} dx - \int_{M} u_{\epsilon} \Delta v dx = 0
\end{equation}
Choosing $v \equiv 1$, we get:

\begin{equation}
\int_{M} \Delta u_{\epsilon} dx = 0
\end{equation}
Hence, we integrate the equation where $u_{\epsilon}$ is known \textit{a priori} to be continuous:

\begin{equation}
\int_{M} \Delta u_{\epsilon} + f(x) + \epsilon u_{\epsilon} dx = 0
\end{equation}

Thus, since $f$ is already chosen to integrate to zero, we get $\int_{M} u_{\epsilon} dx = 0$.

Fix $\epsilon>0$. We will use interior regularity estimates for subsets of Euclidean space to show that the Lipschitz constant of $u_{\epsilon}$ is uniformly bounded. The results in Euclidean space can be used for the compact surface $M$ by using a $C^{\infty}$-atlas to cover the manifold and applying the Euclidean interior estimates in coordinate patches $\Omega_i$ that overlap by a fixed distance $\delta >0$ if necessary. Expressed over small enough coordinate patches in local coordinates, we use estimates that are valid for second-order uniformly elliptic linear PDE with bounded first- and zeroth-order coefficients. We use these estimates, because the Poisson equation expressed in a small enough coordinate patch in local coordinates has first-order terms due to the non-zero Christoffel symbols. However, due to curvature bounds, the coordinate patches can be chosen small enough (satisfying $\text{diam}(\Omega) \leq r_{M}$, where $r_{M}$ is the injectivity radius of the manifold $M$), these Christoffel symbols can be uniformly bounded. Interior Schauder estimates of uniformly elliptic linear PDE in non-divergence form with bounded coefficients in subsets of Euclidean space can be found in Corollary 6.3 of Gilbarg and Trudinger~\cite{gilbargtrudinger}:

\begin{equation}\label{eq:schauder}
d \left\Vert Du_{\epsilon}\right\Vert_{C^{0}(\Omega'_i)} + d^2 \left\Vert D^2u_{\epsilon}\right\Vert_{C^{0}(\Omega'_i)} \leq C \left( \left\Vert u_{\epsilon} \right\Vert_{C^{0}(\Omega_i)} + \left\Vert f \right\Vert_{C^{0}(\Omega_i)} \right) 
\end{equation}
where $d \leq \text{dist}(\Omega ', \partial \Omega)$. In order to avoid the effect of the terms $d$ on the estimate, the result over the compact manifold $M$ will be effected by choosing overlapping coordinate patches $\Omega_i$ that have a uniform constant of overlap $r_{M}>\delta>0$. That is, $d(x,y) > \delta$ for any $x \in M \setminus \Omega_i, y \in M \setminus \Omega_j$ for $i \neq j$. Thus, we will be able to choose $d \geq \delta$. Doing this, we derive an estimate over the whole compact manifold $M$:

\begin{equation}
\delta \left\Vert Du_{\epsilon}\right\Vert_{C^{0}(M)} + \delta^2 \left\Vert D^2u_{\epsilon}\right\Vert_{C^{0}(M)} \leq C \left( \left\Vert u_{\epsilon} \right\Vert_{C^{0}(M)} + \left\Vert f \right\Vert_{C^{0}(M)} \right) 
\end{equation}

The term $\left\Vert u_{\epsilon} \right\Vert_{C^{0}(\Omega_i)}$ on the right-hand side of \eqref{eq:schauder} can be bounded by using the Krylov-Safonov Harnack inequality in Euclidean space, see~\cite{cabre}. Denote $v_{\epsilon} = u_{\epsilon} - \inf_{M} u_{\epsilon}$. Then, for any ball $B_r$ we have:

\begin{equation}
\sup_{B_{r}} v_{\epsilon} \leq C \left( \inf_{B_{r}} v_{\epsilon} + \left\Vert f \right\Vert_{L^2(B_{2r})} \right)
\end{equation}
where $B_{2r}$ denotes the concentric ball with twice the radius of $B_r$. Since we can achieve a finite covering of $M$ by balls $B_{\delta}$, where $\delta<r_M$ as explained above, we now use a chaining argument to get a uniform bound on $\sup_{B_{\delta}} v_{\epsilon}$. Denoting the ball $B^{1}_{\delta}$ where the minimum is obtained, i.e. $v_{\epsilon}(x) = 0$ for some $x \in B^{1}_{\delta}$, we have:

\begin{equation}
\sup_{B^{1}_{\delta}} v_{\epsilon} \leq C \left\Vert f \right\Vert_{L^2(M)}
\end{equation}

Now take a ball $B^{2}_{\delta}$ that overlaps by $\delta$ with $B^{1}_{\delta}$. Then,

\begin{equation}
\sup_{B^{1}_{\delta}} v_{\epsilon} \leq C \left( \inf_{B^{1}_{\delta}} v_{\epsilon} + \left\Vert f \right\Vert_{L^2(M)} \right) \leq C \left( \sup_{B^{1}_{\delta}} v_{\epsilon} + \left\Vert f \right\Vert_{L^2(M)} \right) \leq C' \left\Vert f \right\Vert_{L^2(M)} 
\end{equation}

Repeating this chaining argument aa finite number of times can cover the compact manifold $M$. Thus, we get:

\begin{equation}
\sup_{M} v_{\epsilon} \leq C'' \left\Vert f \right\Vert_{L^2(M)}
\end{equation}
Thus, since $u_{\epsilon}$ is mean-zero, we have therefore that $u_{\epsilon}$ is bounded and thus $\left\Vert u_{\epsilon} \right\Vert_{C^{2,0}(M)} \leq C$ for a universal constant $C$ that does not depend on $\epsilon$. Hence, we know that since the $u_{\epsilon}$ is bounded in its Lipschitz constant, so, in fact, the solutions coincide, i.e. $w_{\epsilon} = u_{\epsilon}$, for $R$ chosen large enough.

By the stability of viscosity solutions~\cite{CIL}, that $\overline{U} = \limsup_{\epsilon} u_{\epsilon}$ is a viscosity subsolution of $- \Delta u + f(x)$. Likewise, we can define $\underline{U} = \liminf u_{\epsilon}$ and it is a viscosity supersolution of $- \Delta u + f(x)$.

Take any sequence $\epsilon_n \rightarrow 0$ such that $u_{\epsilon} \rightarrow U$ uniformly. Then, $U$ is the mean-zero viscosity solution of $- \Delta u + f(x)$. For any sequence $\epsilon_n$, there exists a subsequence $\epsilon_{n_{k}}$ such that $u_{\epsilon_{n_{k}}} \rightarrow U$ uniformly, by Arzela-Ascoli. Since every sequence contains a subsequence that converges to $U$, we conclude that, in fact $u_{\epsilon_{n}} \rightarrow U$ for any sequence.
\end{proof}
\end{thm}
}

\end{document}